\documentclass[ejs]{imsart}

\RequirePackage[OT1]{fontenc}
\RequirePackage{amsthm,amsmath}
\usepackage{color}
\RequirePackage[authoryear]{natbib}
\RequirePackage[colorlinks,citecolor=blue,urlcolor=blue]{hyperref}
\usepackage{epsfig}
\pubyear{2018}
\volume{12}
\issue{0}
\firstpage{2709}
\lastpage{2742}

\startlocaldefs
\numberwithin{equation}{section}
\theoremstyle{plain}
\newtheorem{thm}{Theorem}[section]
\endlocaldefs

\usepackage{amsfonts}
\newtheorem{prop}{Proposition}[section]
\newtheorem{lem}{Lemma}[section]

\newtheorem{rem}{Remark}[section]
\def\ind{\overset{d}{=}}
\def\as{\overset{a.s.}{\to}}
\def\inp{\overset{p}{\to}}
\newcommand{\trans}{^{\mbox{\tiny{T}}}}
\def\pr{\mbox{pr}}

\def\defby{\stackrel{\mbox{\textrm{\tiny def}}}{=}}
\def\tr{\mbox{tr}}
\def\r{\mbox{ER}}
\def\r{\mbox{R}}
\def\mR{\mathbb{R}}
\def\I{\mathbb{I}}

\begin{document}

\begin{frontmatter}
\title{On the dimension effect of regularized linear discriminant analysis} 
\runauthor{Wang and Jiang}
\runtitle{Dimension effect of LDA and RLDA}

\begin{aug}
\author{\fnms{Cheng} \snm{Wang}\thanksref{a1}\ead[label=e1]{chengwang@sjtu.edu.cn}}
\and 
\author{\fnms{Binyan} \snm{Jiang}\thanksref{a2}\ead[label=e2]{by.jiang@polyu.edu.hk}}


\address[a1]{School of Mathematical Sciences, \\
	Shanghai Jiao Tong University, Shanghai, 200240,  China.\\
\printead{e1}}

\address[a2]{Department of Applied Mathematics, \\The Hong Kong Polytechnic University,
	Hung Hom, Kowloon, Hong Kong. \\
	\printead{e2}}

\end{aug}

\begin{abstract}
This paper studies the dimension effect of the linear discriminant analysis (LDA) and the regularized linear discriminant analysis (RLDA) classifiers for large dimensional data where the observation dimension $p$ is of the same order as the sample size $n$. More specifically, built on properties of the Wishart distribution and recent results in random matrix theory, we derive explicit expressions for the asymptotic misclassification errors of LDA and RLDA respectively, from which we gain insights of how dimension affects the performance of classification and in what sense. Motivated by these results, we propose adjusted classifiers by correcting the bias brought by the unequal sample sizes. The bias-corrected LDA and RLDA classifiers are shown to have smaller misclassification rates than LDA and RLDA respectively.  Several interesting examples are discussed in detail and the theoretical results on dimension effect are illustrated via extensive simulation studies.
\end{abstract}

 
\tableofcontents
\end{frontmatter}

\section{Introduction}

Discriminant analysis that aims to allocate objects into one of the predefined classes has been an important topic in statistical learning and data analysis. Modern data is frequently featured by complex structures such as high dimensionality. In the past decades, extensive research has been done to address new challenges of high dimensionality in discriminant analysis.  In particular, due to its simplicity, optimality \citep{anderson2003}, and promising performance in real data analysis \citep{hand2006classifier}, the linear discriminant analysis (LDA) classifier has received more and more attention for classifying data of very large dimension under certain sparsity assumptions \citep{cai2011direct, mai2012direct,fan2012road}. In this paper, we will focus on the general case where no sparsity assumption is assumed and study the effect of dimensionality on the performance of LDA and the regularized LDA \citep{friedman1989regularized}. For simplicity, we shall focus on the two-class classification problem.

The issue of dimension effect, or more specifically, the effect that a diverging dimension $p$ would brought to classical statistical inference \citep{bai2009enhancement}, has been studied in several statistical problems in the literature. \cite{bai1996effect} first studied the dimension effect of Hotelling's $T^2$ test by deriving the explicit power under local alternatives. \cite{bai2009enhancement} and \cite{el2010high} considered the dimension effect of Markowitz's portfolio. For LDA, when the dimension $p$ is assumed to be fixed, many works have been done to study the asymptotic of the misclassification rate by studying the properties of Wishart distribution.  For a diverging $p$, a seminal work by \cite{bickel2004} showed that when $p/n \to \infty$, LDA tends to random guessing.  \cite{shao2011} proved that when $p \sqrt{\log{p}}/\sqrt{n} \to 0$, LDA is consistent in that the error rate of sample LDA tends to the oracle Bayes error rate. Other related works that handle variations of LDA (for example, quadratic discriminant analysis etc.) can also be found in \cite{saranadasa1993asymptotic}, \cite{cheng2004asymptotic} and \cite{li2016two}, among others. On the other hand, when $p$ is large, regularization method is commonly used to reduce the dispersion of the sample covariance matrix; see for example \cite{chen2011regularized} for the study of the regularized Hotelling's $T^2$ test and \cite{ledoit2004honey} for regularized estimation in Markowitz's portfolio. Regularization procedure has also been widely applied in many other statistical analysis in current high dimensional literature; see for example \cite{cai2011constrained},  \cite{buhlmann2013statistical} and \cite{wang2016high}. For discriminant analysis, the regularized LDA (RLDA) was proposed and studied in \cite{friedman1989regularized} and \cite{guo2007regularized}. 
Recently, \cite{zollanvari2013application, zollanvari2015generalized}  studied the misclassification error rate of RLDA with the motivation to estimate this error rate and moreover, they did not provide explicit expressions for the rate. \cite{dobriban2015high} analyzed the predictive risk of ridge regression and regularized discriminant analysis  and derived some explicit formulas. However, the results were based on the random effects hypothesis which is a strong assumption.

Among these existing literatures, the problem of how the high dimensionality would affect the classification accuracy of LDA and RLDA when the observation dimension $p$ is of the same order as the sample size $n$ is still not well understood and in this work we provide a comprehensive analysis of the misclassification error rates under mild conditions. Built on properties of the Wishart distribution, we will first study the dimension effect of LDA (Section 2) under the assumption that $p/n\rightarrow y\in(0,1)$. In the general case where $p/n\rightarrow y\in (0,\infty)$, using recent results in random matrix theory, we provide a systematic study on the dimension effect of the RLDA in Section 3. 
Overall, this paper aims at providing theoretical studies on the dimension effect of LDA and RLDA, from which we gain insights of how the increasing dimension affects the classification accuracy of LDA and RLDA. Following is a summary of the contributions of our work:

\begin{itemize}

\item For LDA, we study the effect of the sample means and the sample covariance matrix respectively. Interestingly, our results show that the dimension effect brought by the sample means would introduce some bias to LDA, which would be further amplified by the effect brought by the sample covariance matrix. Motivated by these observations, a bias correction is proposed to improve the classification accuracy.

\item For RLDA, we derive the closed form of the error rates under general cases. \cite{dobriban2015high} studied a similar question under random effect which is a strong condition. We relax this condition and consider a general case which includes random effect as a special case. Similar to LDA,  a bias corrected RLDA is proposed to reduce the bias brought by the dimension effect that improves the classification performance.

\item To study the dimension effect of RLDA, the key question is to study the asymptotic limits of the moments of a class of random matrices and their random quadratic forms involved of the population means. The former one is studied by \cite{ledoit2011eigenvectors}, \cite{chen2011regularized} and our previous work \cite{wang2015shrinkage}. When the population means satisfy some special structures such as the eigenvector structures in \cite{bai2007asymptotics} or the random effect of \cite{dobriban2015high}, the limits of quadratic forms are the same as the moments. By studying the limits of these moments, \cite{dobriban2015high} derived the dimension effect of RLDA under the random effect assumption. For general cases, the asymptotic of the random quadratic forms are much more challenging.  \cite{bai2011asymptotic} derived results for an identity matrix. \cite{karoui2011geometric} provided a tool to study the limits of the quadratic forms. Based on the methods of \cite{karoui2011geometric}, we develop explicit results for the dimension effect of RLDA under mild conditions.

\end{itemize}

The remainder of this paper is organized as follows. We study the dimension effect of LDA when $p/n\rightarrow y\in(0,1)$ in Section 2. Dimension effect of RLDA is provided in Section 3. To better understand the dimension effect of RLDA and verify our theoretical results, we conduct several interesting examples in Section 4 and simulations in Section 5. All technical details are relegated to the Appendix.

\section{Linear Discriminant Analysis}
 
Let $X$ be a $p$-dimensional normal random vector belonging to class $k$ if $X \sim N(\mu_k,\Sigma),~k=1,2,$ where $\mu_1 \neq \mu_2$ are the population means and $\Sigma$ is the population covariance matrix. For simplicity, in this paper we consider the case where equal weights are used for the two types of misclassification error.
If $\mu_1,\mu_2$ and $\Sigma$ are known, the Bayes' classification rule is 
\begin{align}\label{Bayesrule} 
\mbox{Bayes' rule:}~~\delta_{0}(X)=\I \left\{\Big(X-\frac{\mu_1+\mu_2}{2}\Big) \trans \Sigma^{-1}(\mu_1-\mu_2)>0\right\},
\end{align}
where $\I(\cdot)$ is the indicator function which assigns $X$ to class 1 if and only if $\delta_0(X)=1$. Bayes' rule is known to be optimal in that it has the minimum misclassification rate among all classifiers.  Specifically,  the optimal Bayes error rate is
\begin{align}\label{R0} 
\r_{0}=&\frac{1}{2} \pr \{ \delta_{0}(X)=0|X \sim N(\mu_1,\Sigma)\}+\frac{1}{2} \pr \{ \delta_{0}(X)=1|X \sim N(\mu_2,\Sigma)\}\nonumber\\
=&\Phi\Big(- \frac{\Delta}{2}\Big), 
\end{align}
where $\Delta=\sqrt{(\mu_1-\mu_2) \trans \Sigma^{-1}(\mu_1-\mu_2)}$ and $\Phi(\cdot)$ is the standard normal distribution function.

Let $\{X_{1,j},j=1,\cdots,n_1\}$ and $\{X_{2,j},j=1,\cdots,n_2\}$ be independent and identically distributed random samples from $N_p(\mu_1,\Sigma)$ and $N_p(\mu_2,\Sigma)$, respectively. We can estimate $\mu_1,\mu_2$ and $\Sigma$ by their sample analogs,
\begin{align*}
&\bar{X}_k=\frac{1}{n_k}\sum_{j=1}^{n_k} X_{k,j},~k=1,2,~S_n=\frac{1}{n-2}\sum_{k=1}^2 \sum_{j=1}^{n_k}(X_{k,j}-\bar{X}_k)(X_{k,j}-\bar{X}_k) \trans,
\end{align*}
where $n=n_1+n_2$. Plugging the estimation into Bayes' rule \eqref{Bayesrule} we obtain the classical LDA classifier
\begin{align} \label{fisherLDA}
\delta_{\text{LDA}}(X)=\I\left\{ \left(X-\frac{\bar{X}_1+\bar{X}_2}{2}\right)\trans S_n^{-1}(\bar{X}_1-\bar{X}_2)>0\right\}.
\end{align}
After some simple calculation we can obtain that conditional on the samples the misclassification rate of LDA is
\begin{align*}
\r_{\text{LDA}}=\frac{1}{2}\sum_{j=1}^2 \Phi\left(\frac{(-1)^j(\mu_j-\frac{\bar{X}_1+\bar{X}_2}{2})\trans S_n^{-1}(\bar{X}_1-\bar{X}_2)}{\sqrt{(\bar{X}_1-\bar{X}_2) \trans S_n^{-1} \Sigma S_n^{-1}(\bar{X}_1-\bar{X}_2)}}\right).
\end{align*}
To investigate the dimension effect of the sample means and the sample covariance matrix separately, we will first study two simple cases to gain some insights: 
\begin{itemize}
	\item[(i)] Assuming that $\Sigma$ is known, what is the effect brought by the sample mean?
	\item[(ii)] Assuming that $\mu_1, \mu_2$ are known, what is the effect brought by the sample covariance matrix? 
\end{itemize}
For convenience, we assume that $\Delta$ is a constant throughout the paper.  As pointed out by \cite{cai2011direct}, 
the case where $\Delta\rightarrow 0$ would indicate that no classifier is better than random guessing. Similarly, the case where $\Delta\rightarrow \infty$ would imply that the two classes are so separated that classification is trivial such that even naive Bayes would asymptotically have zero misclassification rate \citep{aoshima2014distance}.
 
\subsection{Effect of the sample mean}
To investigate the effect brought by the sample mean, we first assume $\Sigma$ is known and consider the following linear classifier 
\begin{align*}
\delta_{\text{LDA1}}(X)=\mathbb{I}\left\{\left(X-\frac{\bar{X}_1+\bar{X}_2}{2}\right)\trans \Sigma^{-1}(\bar{X}_1-\bar{X}_2)>0\right\},
\end{align*}
which is obtained by replacing $S_n$ in the LDA classifier \eqref{fisherLDA} by the true $\Sigma$.
The misclassification rate of the above classifier can be computed as: 
\begin{align*}
\r_{\text{LDA1}}  = \frac{1}{2}\sum_{j=1}^2 \Phi\left(\frac{(-1)^j\left(\mu_j-\frac{\bar{X}_1+\bar{X}_2}{2}\right)\trans \Sigma^{-1}(\bar{X}_1-\bar{X}_2)}{\sqrt{(\bar{X}_1-\bar{X}_2) \trans \Sigma^{-1}(\bar{X}_1-\bar{X}_2)}}\right).
\end{align*}
The dimension effect brought by the sample means can then be characterized by the asymptotic misclassification rate given in the following theorem.
\begin{thm} \label{thm1} Assuming $p/n_j \to y_j \in (0,\infty),~j=1,2$, we have,
	\begin{align}\label{R1}
\r_{\text{LDA1}} \inp \frac{1}{2} \Phi\Bigg(- \frac{\Delta^2+y_2-y_1}{2\sqrt{\Delta^2+y_1+y_2}}\Bigg)+\frac{1}{2}\Phi\Bigg(- \frac{\Delta^2+y_1-y_2}{2\sqrt{\Delta^2+y_1+y_2}}\Bigg).
	\end{align}
\end{thm}
By comparing \eqref{R1} with the optimal Bayes error \eqref{R0}, we can see that the dimension effect brought 
by the sample means depends on $y_1, y_2$ only. From Theorem \ref{thm1} we know that even if $\Sigma$ is known, when $y_1$ or $y_2$ is very large which indicates that the dimension far more exceeds the sample size, the misclassification rate would tend to $50 \%$.  The conclusion is consistent with the one in the seminal paper \cite{bickel2004} which studied the dimension effect from a minimax point of view under the setting $p/n \to \infty$. 
\subsection{Effect of the sample covariance matrix}
Similarly, we study the effect brought by the sample covariance matrix. Assuming the population means $\mu_1,\mu_2$ are known, we consider the classifier
\begin{align*}
\delta_{\text{LDA2}}(X)=\mathbb{I}\left\{\left(X-\frac{\mu_1+\mu_2}{2}\right) \trans  S_n^{-1} (\mu_1-\mu_2)>0\right\}.
\end{align*}
The corresponding misclassification rate is 
\begin{align*}
\r_{\text{LDA2}}=\Phi\left(\frac{-(\mu_1-\mu_2)\trans S_n^{-1} (\mu_1-\mu_2) }{2\sqrt{(\mu_1-\mu_2)\trans S_n^{-1}\Sigma S_n^{-1} (\mu_1-\mu_2)}}\right).
\end{align*}
\begin{thm} \label{thm2}
	Assuming $p/(n_1+n_2) \to y \in (0,1)$, we have,
	\begin{align}\label{R2}
\r_{\text{LDA2}} \inp  \Phi\left(- \frac{\Delta}{2} \sqrt{1-y}\right).
	\end{align} 
\end{thm}
Comparing with the optimal Bayes error \eqref{R0} we can see that the term $\sqrt{1-y}$ in \eqref{R2} is exactly the price we pay by using the sample covariance matrix in LDA. Since the results only depend on Wishart distribution, it would be possible to built up more precise results under weaker conditions such as $n>p+2$ and $p/n \to y \in (0,1]$ \citep{jiang2013central} and we leave it to future work. By the properties of Wishart distribution \citep{cook2011mean}, we have
\begin{align*}
E S_n^{-1}=\frac{n-2}{n-p-3}\Sigma^{-1},~E S_n^{-1} \Sigma S_n^{-1}=\frac{(n-2)^2(n-3)}{(n-p-2)(n-p-3)(n-p-5)}\Sigma^{-1}.
\end{align*}
Therefore, the term $\sqrt{1-y}$ is actually introduced by the inverse Wishart distribution. We note  $\sqrt{1-y}$ 
 also arises in the dimension effect for Hotelling's test  \citep{pan2011central} or the Markowitz portfolio \citep{bai2009enhancement, el2010high} which both involve the inverse of the sample covariance matrix.

\subsection{Dimension effect of LDA}
Theorems \ref{thm1} and \ref{thm2} show the explicit results for the sample means and the sample covariance from which we can see the explicit price we pay for the estimation. Now, we present the result for LDA as follows.

\begin{thm} \label{thm3}
Assuming $p/n_j \to y_j \in (0,\infty),~j=1,2$ and $y \defby y_1y_2/(y_1+y_2)<1$, we have
	\begin{align}\label{MR:lda}
	\r_{\text{LDA}}\inp \frac{1}{2}\sum_{j=1}^2 \Phi\left(- \frac{\Delta^2+(-1)^j(y_1-y_2)}{2\sqrt{\Delta^2+y_1+y_2}} \sqrt{1-y}\right).
	\end{align}	
\end{thm}
Note that the term on the right hand side of \eqref{MR:lda} correspond to the misclassification rate within class 1 and class 2 respectively.  LDA has different classification performance on the two classes when $y_1\neq y_2$. From the proof of Theorem \ref{thm1} we can see that this is due to the estimation bias of intercept part in LDA. On the other hand, from \eqref{R2} we learn that the effect from the sample covariance matrix would result in a multiplication factor $\sqrt{1-y}$. Interestingly, we show that when both sample means and sample covariance matrix are used as in LDA \eqref{fisherLDA}, the bias introduced by the sample means was further amplified by the multiplication factor $\sqrt{1-y}$.

To investigate the bias of LDA, we consider the classifier,
\begin{align}\label{general}
\I\left\{  X \trans  S_n^{-1}(\bar{X}_1-\bar{X}_2)+\alpha>0\right\}.
\end{align} 
By Proposition 2 of \cite{mai2012direct}, given the classification direction $S_n^{-1}(\bar{X}_1-\bar{X}_2)$ in \eqref{general}, the optimal intercept corresponding to minimum misclassification rate is given as,
\begin{align*}
\alpha_0=-\frac{1}{2} (\mu_1+\mu_2) \trans S_n^{-1}(\bar{X}_1-\bar{X}_2).
\end{align*}
One the other hand, note that \eqref{general} would reduce to the LDA classifier \eqref{fisherLDA} if the intercept $\alpha$ is set to be  
$
-\frac{1}{2} (\bar{X}_1+\bar{X}_2) \trans S_n^{-1}(\bar{X}_1-\bar{X}_2).
$
However, direct calculation shows that, 
\begin{align*}
-\frac{1}{2}  E (\bar{X}_1+\bar{X}_2) \trans S_n^{-1}(\bar{X}_1-\bar{X}_2)=E \alpha_0+\frac{n-2}{n-p-3}\Big(\frac{p}{2 n_2}-\frac{p}{2 n_1}\Big),
\end{align*}
which indicates that the expected bias between the intercept of LDA and the optimal intercept $\alpha_0$ is exactly $\frac{n-2}{n-p-3}\big(\frac{p}{2 n_2}-\frac{p}{2 n_1}\big)$.   
We thus can make a bias correction to LDA as follow,
\begin{align} \label{clda}
\delta_{\text{LDA}}^c(X)=\mathbb{I}\left\{\Big(X-\frac{\bar{X}_1+\bar{X}_2}{2}\Big)\trans S_n^{-1}(\bar{X}_1-\bar{X}_2)+\frac{n-2}{n-p-3}\Big(\frac{p}{2 n_1}-\frac{p}{2 n_2}\Big)>0  \right\}.
\end{align}
Let $\r_{\text{LDA}}^c$ be the misclassification rate of the bias-corrected classifier \eqref{clda}. The following proposition gives the asymptotic misclassification rate of the above bias-corrected classifier.  
\begin{prop} \label{prop0}
	Under the conditions of Theorem \ref{thm3},  
	\begin{align}\label{MR:clda}
\r_{\text{LDA}}^c \inp \Phi\left(- \frac{\Delta^2}{2\sqrt{\Delta^2+y_1+y_2}} \sqrt{1-y}\right).
	\end{align}
\end{prop}
By noticing that $\Phi(x)$ is strictly convex in $x\in(-\infty, 0)$, we immediately have that the asymptotic misclassification rate of the bias corrected LDA given in \eqref{MR:clda} is smaller than the asymptotic misclassification rate of LDA given in \eqref{MR:lda}. In addition, considering the linear classifier with optimal intercept, 
\begin{align*}
\I\left\{  X \trans  S_n^{-1}(\bar{X}_1-\bar{X}_2)-\frac{1}{2} (\mu_1+\mu_2) \trans S_n^{-1}(\bar{X}_1-\bar{X}_2)>0\right\}.
\end{align*} 
It can be shown that the result given in \eqref{MR:clda} would also hold for this oracle classifier, indicating that our bias correction procedure does eliminate the bias introduced by the unequal sample sizes. Finally, we remark that bias issues in binary classification have received much attention in early literatures; see for example \cite{chan2009scale} and  \cite{huang2010bias} and the references therein. For LDA, the bias-corrected classifier \eqref{clda} happens to be identical to formula (8) of \cite{moran1979closer}. Our motivation and interpretation provide a different view on \eqref{clda} and theoretical justifications provided in Proposition \ref{prop0} are also new.  

In literature, to handle high dimension data, 
naive Bayes \citep{dudoit2002comparison, bickel2004} and sparse LDA \citep{cai2011direct, mai2012direct, fan2012road} were proposed. For naive Bayes, the correlation between the covariates is  ignored and only the diagonal elements of the sample covariance matrix are used.  Sparse LDA methods assumed that the discriminant direction $\beta=\Sigma^{-1}(\mu_1-\mu_2)$ is sparse in that only several elements in $\beta$ are non-zero and others are either zero or close to zero. In the following remarks we provide some comparison of the LDA rule to the naive Bayes rule and sparse LDA approaches.

\begin{rem} \label{nbayes}
For simplicity, we consider the homogeneous case where the diagonal elements of $\Sigma$ equal to each other and assume that
$n_1=n_2$. Consider 
the following naive Bayes classification rule,
\begin{align*}
\delta_{\text{NB}}(X)=\mathbb{I}\left\{\left(X-\frac{\bar{X}_1+\bar{X}_2}{2}\right)\trans (\bar{X}_1-\bar{X}_2)>0\right\}.
\end{align*}
 The conditional misclassification rate of $\delta_{\text{NB}}(X)$ is
\begin{align*}
\r_{\text{NB}}  = \frac{1}{2}\sum_{j=1}^2 \Phi\left(\frac{(-1)^j\left(\mu_j-\frac{\bar{X}_1+\bar{X}_2}{2}\right)\trans (\bar{X}_1-\bar{X}_2)}{\sqrt{(\bar{X}_1-\bar{X}_2) \trans \Sigma (\bar{X}_1-\bar{X}_2)}}\right).
\end{align*}
By a similar analysis as our proofs for Theorem \ref{thm3}, we have 
\begin{align*}
\r_{\text{NB}} = \Phi\left(-\frac{(\mu_1-\mu_2)\trans (\mu_1-\mu_2)+o_p(1) }{2 \sqrt{(\mu_1-\mu_2) \trans \Sigma (\mu_1-\mu_2)+\frac{4}{n} \tr (\Sigma^2)}+o_p(1) }\right).
\end{align*}
We can see that the performance of naive Bayes depends on the structure of $\Sigma$. When $\Sigma=\sigma^2I_{p}$ for some constant $\sigma^2>0$, comparing with \eqref{MR:lda}, the asymptotic error rate of the naive Bayes rule is the same as that of LDA except the term $\sqrt{1-y}$, indicating that naive Bayes is always better than LDA. When the off diagonal elements of $\Sigma$ are nonzero, we would need to pay some price for wrongly ignoring the correlations. A simulation study will be conducted in Section 5, providing further numerical comparison between naive Bayes and LDA. 
\end{rem}

\begin{rem}\label{slda}
For sparse LDA, we take the linear programming discriminant (LPD) rule as an example which was proposed by \cite{cai2011direct}.  The LPD rule estimates the discriminant direction by solving the following optimization problem,
\begin{align*}
\hat{\beta}= \arg \min |\beta|_1,~\mbox{subject~to~} |S_n \beta-(\bar{X}_1-\bar{X}_2)|_\infty \leq t_n,
\end{align*} 
where $t_n$ is a tuning parameter. Under regular conditions, the misclassification rate $R_{\text{LPD}}$ of LPD \citep[Theorem 3]{cai2011direct} satisfies,
\begin{align*}
\frac{R_{\text{LPD}}}{\Phi(-\Delta/2)}-1=O_p\left(\Delta^2 s \sqrt{\frac{\log{p}}{n}}\right).
\end{align*}
where $s$ is the number of non-zero elements of $\Sigma^{-1}(\mu_1-\mu_2)$. Unlike LDA, the misclassification rate of LPD not only depends on $p$ and $n$, but also hinge on the size of $s$. When $s$ is small, LPD will achieve a good performance and even close to the Bayes rule while the dimension $p$ is allowed to grow exponentially in $n$. However, when the sparsity assumption is violated, it is expected that  LPD will perform poor if $s$ is large.  Further numerical comparison will be provided in Section 5.
\end{rem}

\section{Regularized Linear Discriminant Analysis}
RLDA \citep{friedman1989regularized, guo2007regularized} was proposed to reduce the dispersion of eigenvalues of the sample covariance matrix when $p$ is large and overcome the singularity issue when $p>n-2$. For a given tuning parameter $\lambda>0$, the RLDA classifier is given as,
\begin{align*}
\delta_{\text{RLDA}}(X)=\mathbb{I}\left\{\left(X-\frac{\bar{X}_1+\bar{X}_2}{2}\right)\trans (S_n+\lambda I_p)^{-1}(\bar{X}_1-\bar{X}_2)>0\right\}.
\end{align*}
When $S_n$ is invertible, RLDA would reduce to LDA if $\lambda$ is set to be zero. With one more degree of freedom to adjust the estimation of the covariance matrix, it is well known that RLDA generally poses better performance than LDA, providing that $\lambda$ is chosen properly. The misclassification rate of RLDA  is
\begin{align*}
\r_{\text{RLDA}}(\lambda)=\frac{1}{2}\sum_{j=1}^2 \Phi\left(\frac{(-1)^j\left(\mu_j-\frac{\bar{X}_1+\bar{X}_2}{2}\right)\trans (S_n+\lambda I_p)^{-1}(\bar{X}_1-\bar{X}_2) }{\sqrt{(\bar{X}_1-\bar{X}_2)\trans (S_n+\lambda I_p)^{-1}\Sigma (S_n+\lambda I_p)^{-1} (\bar{X}_1-\bar{X}_2)}}\right).
\end{align*}
To explore the dimension effect on RLDA, we need to study the asymptotic properties of $\r_{\text{RLDA}}(\lambda)$. Denote
\begin{align}\label{Bn}
\mu=\Sigma^{-\frac{1}{2}} (\mu_1-\mu_2),~B_n(\lambda)=\Sigma^{\frac{1}{2}} (S_n+\lambda I)^{-1} \Sigma^{\frac{1}{2}}.
\end{align}
By noticing that
\begin{align*}
\bar{X}_1\ind \frac{1}{\sqrt{n_1}} \Sigma^{\frac{1}{2}}Y_1+\mu_1,~\bar{X}_2 \ind \frac{1}{\sqrt{n_2}} \Sigma^{\frac{1}{2}}Y_2+\mu_2,
\end{align*}
where $Y_1,Y_2 \sim N(0,I_p)$ and $Y_1,Y_2, S_n$ are independent, we can see that 
the asymptotic properties of $\r_{\text{RLDA}}(\lambda)$ would rely on the asymptotic properties of
\begin{align} \label{m0}
\frac{1}{p} Y_1 \trans \tr B_n(\lambda) Y_1 ,~~\frac{1}{p} Y_1 \trans \tr B^2_n(\lambda) Y_1,
\end{align} 
and the random quadratic forms 
\begin{align}  \label{q1}
\mu \trans B_n(\lambda) \mu ,~~ \mu \trans B^2_n(\lambda) \mu.
\end{align}

Under mild conditions, it can be shown that
\begin{align} \label{m1}
\frac{1}{p} Y_1 \trans \tr B^k_n(\lambda) Y_1-\frac{1}{p} \tr B^k_n(\lambda)=o_p(1),~k=1,2.
\end{align}
Thus, finding the limits of the statistics \eqref{m0} can be reduced to finding the moments of $\tr B^k_n(\lambda) /p$. Related studies on the moments of $\tr B^k_n(\lambda) /p$ can be found in \cite{ledoit2011eigenvectors}, \cite{chen2011regularized} and \cite{wang2015shrinkage}. Unlike the statistics \eqref{m0} whose limits only depend on the eigenvalues of $B_n(\lambda)$, the asymptotic of  the quadratic forms \eqref{q1} generally depend on the entry-wise properties of $B_n(\lambda)$ \citep{karoui2011geometric}. Under special case such as $\mu$  or $\mu_1-\mu_2$ itself is generated by i.i.d. entries, the quadratic forms can be expressed by the trace of a random matrix. For example,  \cite{dobriban2015high} studied RLDA by setting $\mu_1-\mu_2$ is i.i.d. where the statistics \eqref{m1} can be approximated by $\tr (S_n+\lambda I_p)^{-k}/p,~k=1,2$, respectively. Here, we study the problem under more general settings where no structure on $\mu_1, \mu_2$ is assumed. 
Before proceed to the theoretical results, we introduce some technical assumptions.  
\begin{enumerate}
	\item[(C1):] $p/n_j \to y_j \in (0,\infty),~j=1,2$ and  denote $y=y_1y_2/(y_1+y_2)$.
	\item[(C2):] The eigenvalues of $\Sigma$ are uniformly bounded and the empirical spectral distribution $F^{\Sigma}$ converges to a nonrandom distribution function $H$ as $p\rightarrow \infty$.
	\item [(C3):] For $t \geq 0$,
	\begin{align*}
	& \Delta^{-2}\mu \trans (I_p+t\Sigma^{-1})^{-1} \mu\to h_1(t),  \\
	&  \Delta^{-2}\mu \trans (I_p+t\Sigma^{-1})^{-2} \mu  \to h_2(t)  .
	\end{align*}
\end{enumerate}
We remark that (C1) and (C2) are common conditions in random matrix theory and (C3) is a technical assumption to express the explicit limits of the quadratic forms \eqref{q1}.  For a given $\lambda$, (C3) can be relaxed to some specific $t$ which is a function of $\lambda$ and here for brevity we make the assumption for general $t$. 
   Under the conditions (C1) and (C2), for any $\lambda>0$, it can be shown that,
\begin{align*}
\frac{1}{p} \tr\{(S_n+ \lambda I_p)^{-1}\} &\overset{a.s.}\longrightarrow m_0(-\lambda), \\
\frac{1}{p} \tr\{(S_n+ \lambda I_p)^{-2}\} &\overset{a.s.}\longrightarrow m'_0(-\lambda)=\frac{d m_0(z)}{d z}\Big|_{z=-\lambda},
\end{align*}
where $m_0(-\lambda)$ is the unique solution of the Mar{\v{c}}enko-Pastur equation \citep{marvcenko1967distribution, el2008spectrum},
\begin{align*}
m(-\lambda)=\int \frac{d H(t)}{t(1-y+y \lambda m(-\lambda))+\lambda},
\end{align*}
under the condition $1-y+y \lambda m(-\lambda)\geq 0$. More details can be found in \cite{wang2015shrinkage}.

 Following is our main results.

\subsection{Dimension effect of RLDA}
\begin{thm} \label{thm4}
Under the conditions (C1)-(C3), for any $\lambda>0$,
	\begin{align}\label{MR:rlda}
	\r_{\text{RLDA}}(\lambda) \inp \frac{1}{2}\ \sum_{j=1}^2  \Phi\left(-\frac{H_1(\lambda) \Delta^2 +(-1)^j(y_1-y_2) R_1(\lambda)}{2 \sqrt{H_2(\lambda) \Delta^2 +(y_1+y_2) R_2(\lambda)}}\right),
	\end{align}	
	where 
\begin{align*}
R_1(\lambda) =& \frac{1-\lambda m_0(-\lambda)}{1-y[1-\lambda m_0(-\lambda)]},\\
R_2(\lambda) =& \frac{1-\lambda m_0(-\lambda)}{\{1-y[1-\lambda
	m_0(-\lambda)]\}^3}-\frac{\lambda m_0(-\lambda)-\lambda^2
	m'_0(-\lambda)}{\{1-y[1-\lambda m_0(-\lambda)]\}^4},\\
H_1(\lambda)=& \frac{1}{1-y[1-\lambda m_0(-\lambda)]}h_1\Bigg(\frac{\lambda}{1-y[1-\lambda m_0(-\lambda)]}\Bigg),\\
H_2(\lambda)=&\{(1+yR_1(\lambda))^2+y R_2(\lambda)\}h_2\left(\frac{\lambda}{1-y[1-\lambda m_0(-\lambda)]}\right).
\end{align*}
\end{thm}
Theorem \ref{thm4} provides an explicit expression for the asymptotic misclassification rate of RLDA, which characterizes the effect of the data dimension in terms of $y_1, y_2$ and the functions on the regularization parameter $\lambda$ defined above. 
\begin{rem}
From the proof of Theorem \ref{thm4}, we note that when $y<1$ and $\lambda=0$, \eqref{MR:rlda} reduces to \eqref{MR:lda}, implying that Theorem \ref{thm3} could be regarded as a special case of Theorem \ref{thm4}. We remark that the proof of Theorem \ref{thm3} is built on properties of normal and Wishart distributions only while Theorem \ref{thm4} is built on the more complicated random matrix theory. 
\end{rem}

By noticing that if 	
\begin{align} \label{con1}
h_1(t)=\int \frac{x}{x+t} dH(x),~h_2(t)=\int \frac{x^2}{(x+t)^2} dH(x),	
\end{align} 
we would have $H_1(\lambda)=R_1(\lambda), H_2(\lambda)=R_2(\lambda)$ and hence the asymptotic misclassification rate \eqref{MR:rlda} can be further simplified. This is summarized in the following proposition.

\begin{prop} \label{prop2}
	Under the conditions of Theorem \ref{thm4} and assuming that \eqref{con1} holds,
	we have
	\begin{align*}
	\r_{\text{RLDA}}(\lambda)\inp \frac{1}{2}\ \sum_{j=1}^2  \Phi\left(-\frac{ \Delta^2+(-1)^j(y_1-y_2)}{2 \sqrt{\Delta^2 +y_1+y_2} }\frac{R_1(\lambda)}{\sqrt{R_2(\lambda)}}\right).
	\end{align*}	
\end{prop}

Note that condition \eqref{con1} would hold when $\mu_1-\mu_2$ and $\Sigma$ satisfy the following structures,
\begin{align*}
\mu \trans (I_p+t\Sigma^{-1})^{-k} \mu=\|\mu\|^2 \frac{1}{p} \tr (I_p+t\Sigma^{-1})^{-k},~k=1,2,
\end{align*} 
A trivial example satisfies the above structure is $\Sigma=I_p$ and an interesting isotropic example will also be discussed in Section 4. 
\subsection{Bias correction for RLDA}
To reduce the bias brought by the unequal sample sizes, 
following the same idea of bias correction for LDA as in Section 2.3, we consider the following class of classifiers,
\begin{align*} 
\I\left\{  X \trans  (S_n+\lambda I_p)^{-1}(\bar{X}_1-\bar{X}_2)+\alpha>0\right\}.
\end{align*} 
Given the discriminant direction $(S_n+\lambda I_p)^{-1}(\bar{X}_1-\bar{X}_2)$, the optimal intercept corresponding to minimum misclassification rate is given as,
\begin{align*}
\alpha_r=-\frac{1}{2} (\mu_1+\mu_2) \trans(S_n+\lambda I_p)^{-1}(\bar{X}_1-\bar{X}_2),
\end{align*}
while for RLDA the intercept is $-\frac{1}{2} (\bar{X}_1+\bar{X}_2) \trans (S_n+\lambda I_p)^{-1}(\bar{X}_1-\bar{X}_2).$  
On the other hand, from the proof of Theorem \ref{thm4} we have 
\begin{align*} 
&-\frac{1}{2}  E (\bar{X}_1+\bar{X}_2) \trans (S_n+\lambda I_p)^{-1}(\bar{X}_1-\bar{X}_2)\\
=&E \alpha_r+E\frac{1}{2n_1}Y_1\trans B_n(\lambda)Y_1-E\frac{1}{2n_2}Y_2\trans B_n(\lambda)Y_2 \\
=&E \alpha_r- \left(\frac{p}{2n_1}-\frac{p}{2n_2}\right) \frac{E \tr(B_n(\lambda))}{p}.
\end{align*}  
Therefore, we can correct the bias as follow,
\begin{align*} 
\left(X-\frac{\bar{X}_1+\bar{X}_2}{2}\right)\trans (S_n+\lambda I_p)^{-1}(\bar{X}_1-\bar{X}_2)+\left(\frac{p}{2n_1}-\frac{p}{2n_2}\right) \frac{E \tr(B_n(\lambda))}{p}>0.
\end{align*}
The matrix $B_n(\lambda)$ depends on the true population covariance matrix $\Sigma$ and the regularized term $(S_n +\lambda I_p)^{-1}$. It is complicated to derive the explicit results for $E \tr(B_n(\lambda))$. By noticing that $\tr(B_n(\lambda))/p$ converges to $R_1(\lambda)$ almost surely \citep{wang2015shrinkage}, we consider estimating $ E\tr(B_n(\lambda))/p$ by the consistent estimator given in \citep{wang2015shrinkage} and propose the following bias-corrected RLDA classifier,
\begin{align}\label{bcrlda}
\delta^c_{\text{RLDA}}(X)=\mathbb{I}&\bigg\{\left(X-\frac{\bar{X}_1+\bar{X}_2}{2}\right)\trans (S_n+\lambda I_p)^{-1}(\bar{X}_1-\bar{X}_2) \nonumber \\
&+\left(\frac{p}{2n_1}-\frac{p}{2n_2}\right)\frac{1-\frac{1}{p}tr(\frac{1}{\lambda}S_n+I_p)^{-1}}{1-\frac{p}{n-2}+\frac{1}{n-2}tr(\frac{1}{\lambda}S_n+I_p)^{-1}}
>0\bigg\}.
\end{align}
Let $\r^c_{\text{RLDA}}(\lambda)$ be the misclassification rate of the bias-corrected classifier\eqref{bcrlda}. 
\begin{prop} \label{prop1}
Under the conditions of Theorem \ref{thm4},	
		\begin{align}\label{MR:bcrlda}
		\r^c_{\text{RLDA}}(\lambda)\inp \Phi\left(-\frac{H_1(\lambda) \Delta^2}{2 \sqrt{H_2(\lambda) \Delta^2 +(y_1+y_2) R_2(\lambda)}}\right).
		\end{align}	
\end{prop}
Again by the convexity of $\Phi(x)$ in $x\in(-\infty, 0)$, we conclude that the asymptotic misclassification rate of the bias corrected RLDA given in \eqref{MR:bcrlda} is smaller than the asymptotic misclassification rate of RLDA given in \eqref{MR:rlda}. Similarly to the bias correction for LDA, we can claim the result of \eqref{MR:bcrlda} is also the asymptotic error rate of RLDA with optimal intercept.

\subsection{Selection of $\lambda$} \label{slambda}
From the perspective of applications, the RLDA can avoid the singularity problem of $S_n$ when $p>n$ and provides a better performance than LDA even if $p<n$.  However, the performance of RLDA relies on the choice of the regularization parameter $\lambda$. Empirically, we can use re-sampling procedures such as cross validation to choose $\lambda$.  However, there is little work on the theoretical analysis of these re-sampling procedures, and cross validation can be computationally expensive when both $p$ and $n$ are large. In this section, we discuss the effect of $\lambda$ on the misclassification error rate and provide a direct estimation of the optimal $\lambda$ for the special cases. 

By Proposition \ref{prop1}, the optimal $\lambda$ with minimum error rate is 
\begin{align}
\lambda_0 \in \arg \max \frac{H_1(\lambda) \Delta^2}{\sqrt{H_2(\lambda) \Delta^2 +(y_1+y_2) R_2(\lambda)}}.
\end{align}
From the definitions of $H_1(\lambda)$ and  $H_2(\lambda)$, these parameters depend on the true structures of $\Sigma$ and $\mu$ and it is not easy to construct estimations for them. However, under the assumption of Proposition \ref{prop2}, the optimal $\lambda$ has a simple form
\begin{align*}
\lambda_0 \in \arg \max \frac{R^2_1(\lambda)}{R_2(\lambda)}.
\end{align*}
This motivates us to select $\lambda_0$ by maximizing an estimate of ${R^2_1(\lambda)}/{R_2(\lambda)}$. As discussed before, $R_1(\lambda)$ and $R_2(\lambda)$ only depends on the eigenvalues of $\Sigma$ and could be estimated by the eigenvalues of $S_n$. Considering a shrinkage estimation $\alpha (S_n+\lambda I_p)^{-1}$ \citep{kubokawa2008estimation, wang2015shrinkage} for the precision matrix $\Sigma^{-1}$, we have 
\begin{align*}
\min_{\alpha} \frac{1}{p}\tr(\alpha (S_n+\lambda I_p)^{-1} \Sigma-I_p)^2 \inp 1-\frac{R^2_1(\lambda)}{R_2(\lambda)}.
\end{align*}
Therefore, similar to \cite{wang2015shrinkage}, we estimate the tuning parameter by 
\begin{align*}
\hat{\lambda}_0= \arg \min \left\{1-\frac{(\hat{R}_1(\lambda))^2}{\hat{R}_2(\lambda)}\right\}, 
\end{align*}
where
\begin{align*}
&\hat{y}=\frac{p}{n-2}, ~a_1(\lambda)=1-\frac{1}{p}tr\Big(\frac{1}{\lambda}S_n+I_p\Big)^{-1},\\ &a_2(\lambda)=\frac{1}{p}tr\Big(\frac{1}{\lambda}S_n+I_p\Big)^{-1}-\frac{1}{p}tr\Big(\frac{1}{\lambda}S_n+I_p\Big)^{-2},\\
&\hat{R}_1(\lambda)=\frac{a_1(\lambda)}{1-\hat{y}a_1(\lambda)},~\hat{R}_2(\lambda)=\frac{a_1(\lambda)}{(1-\hat{y}a_1(\lambda))^3}-\frac{a_2(\lambda)}{(1-\hat{y}a_1(\lambda))^4}. 
\end{align*}
 
Following our previous work \cite{wang2015shrinkage}, it could be possible to establish some theoretical analysis of the estimation $\hat{\lambda}_0$ under the conditions given in the formula \eqref{con1}. For general $\Sigma,\mu_1,\mu_2$, theoretically it is hard to conduct such an effective estimation. Therefore, we skip the  theoretical analysis of the proposed tuning parameter estimation and  conduct simulations in Section 5, to compare the $\hat{\lambda}_0$ proposed above with the one obtained by cross validations.

\section{Examples}
In this section, we use several examples to illustrate the theoretical results obtained in Section 3. For simplicity, in all numerical studies, $n_1$ and $n_2$ are set to be equal and hence the bias-corrected RLDA is identical to RLDA. From Theorem \ref{thm4}, given the ratios $y_1,y_2$, we know the misclassification error rate depends on $R_k(\lambda),k=1,2$ and  $H_k(\lambda),k=1,2$ where $R_k(\lambda)$ is determined by the population covariance $\Sigma$ and $H_k(\lambda)$ involves the structure between $\mu_1-\mu_2$ and $\Sigma$. Write the spectral decomposition of $\Sigma$ as
\begin{align*}
\Sigma=\sum_{i=1}^p \lambda_i v_i v_i \trans,
\end{align*}
where $\lambda_1 \geq \cdots \geq \lambda_p>0$ are the eigenvalues of $\Sigma$ and $v_1,\cdots,v_p$ are the corresponding eigenvectors. To be specific, $R_k(\lambda)$ is determined by the eigenvalues $\lambda_1,\cdots,\lambda_p$ while $H_k(\lambda)$ involves the eigenvalues  as well as the inner products $v_k \trans (\mu_1-\mu_2),k=1,\cdots,p$.  Using  the spectral decomposition of $\Sigma$, we have
	\begin{align*}
& (\mu_1-\mu_2)\trans \Sigma^{-1/2}(I_p+t\Sigma^{-1})^{-1}  \Sigma^{-1/2} (\mu_1-\mu_2)=\sum_{i=1}^p \frac{1}{\lambda_i+t} (v_i \trans (\mu_1-\mu_2))^2,  \\
& (\mu_1-\mu_2)\trans \Sigma^{-1/2}(I_p+t\Sigma^{-1})^{-2}  \Sigma^{-1/2} (\mu_1-\mu_2)=\sum_{i=1}^p \frac{\lambda_i}{(\lambda_i+t)^2} (v_i \trans (\mu_1-\mu_2))^2.
\end{align*}
In what following, we consider three cases.  Firstly, we consider $\Sigma=\sigma^2 I_p$ where all the eigenvalues are identical, and the results would depends on $\sum_{i=1}^p  (v_i \trans (\mu_1-\mu_2))^2=(\mu_1-\mu_2) \trans (\mu_1-\mu_2)$ only. Secondly, we consider three isotropic cases where $\mu_1-\mu_2$, $ \Sigma^{-1/2}(\mu_1-\mu_2)$ or $ \Sigma^{-1}(\mu_1-\mu_2)$ has nonzero projections on all the eigenvectors of $\Sigma$. Consequently, we would have $v_i \trans (\mu_1-\mu_2)\neq 0$ for all $i$. Lastly, we consider the case where $\mu_1-\mu_2$ is parallel to one of the eigenvectors of $\Sigma$, and hence only one of the $v_i \trans (\mu_1-\mu_2),~i=1,\cdots,p$ is nonzero.

\subsection{$\Sigma=\sigma^2 I_p$}
\begin{figure}[!htbp]
	\centerline{
		\begin{tabular}{ccc}				
			\psfig{figure=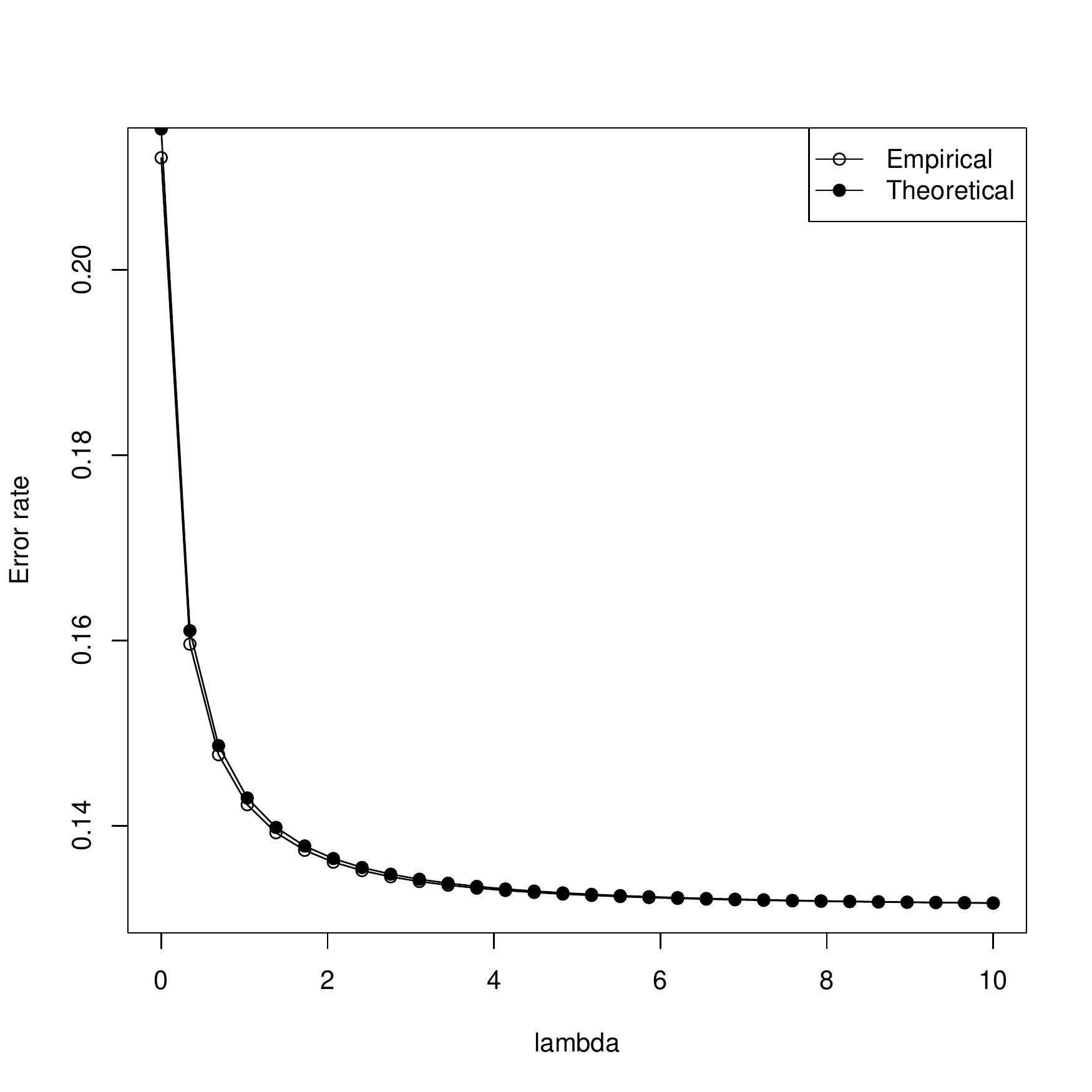,width=1.6 in,height=1.7 in,angle=0} &
			\psfig{figure=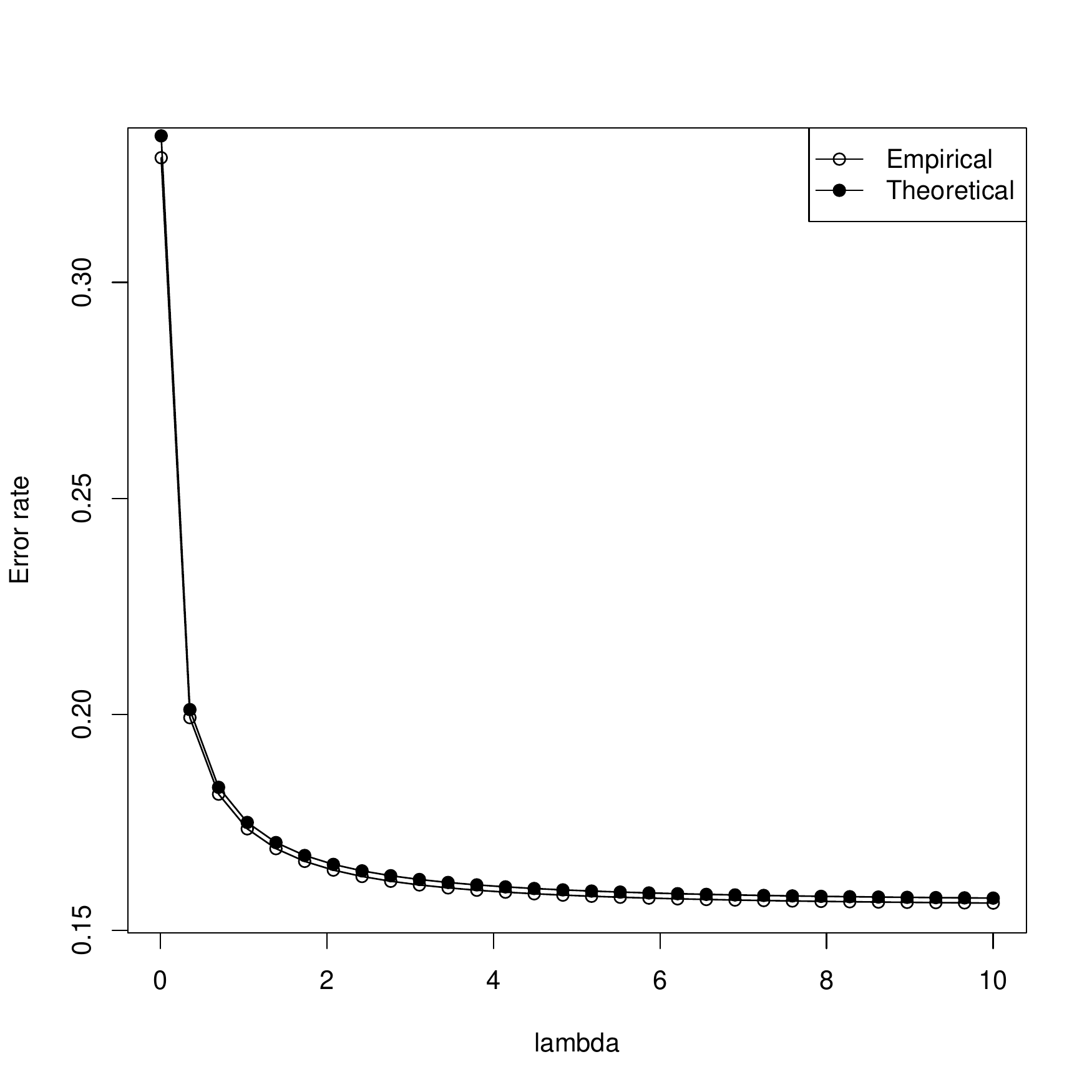,width=1.6 in,height=1.7 in,angle=0} &
			\psfig{figure=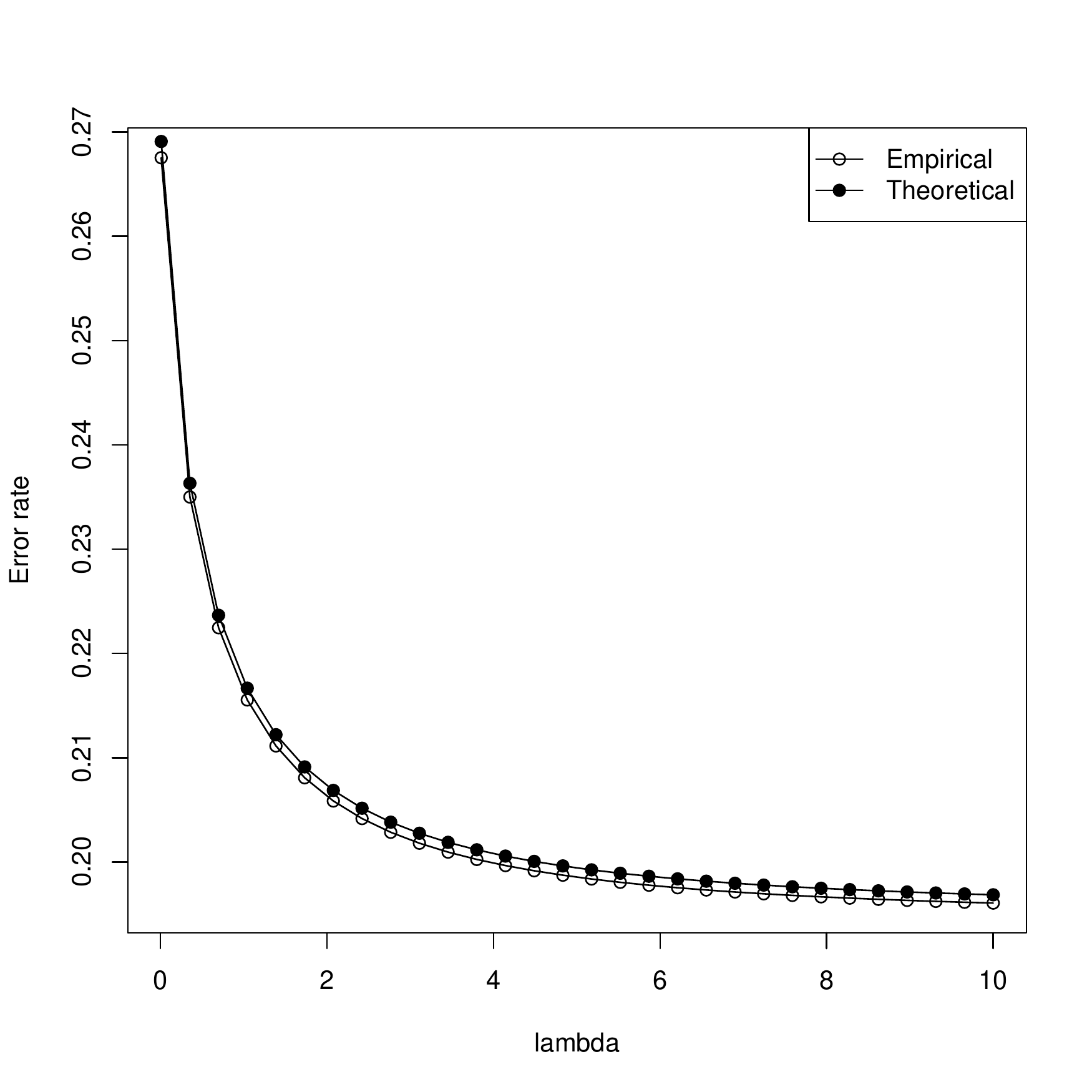,width=1.6 in,height=1.7 in,angle=0}\\
			$p=100$ & $p=200$ & $p=400$
		\end{tabular}
	}
	\caption{Simulations for independent cases where $\Sigma=I_p,~n_1=n_2=100$, $\mu_2=(0,\cdots,0)$ and $\mu_1=(2.56,0,\cdots,0)$. The Bayes misclassification error rate is $10 \% $.}
	\label{fig1}
\end{figure}

We first study the case where $\Sigma$ has a simple form $\Sigma=\sigma^2 I_p$.
It can be easily seen that (C2) and (C3) hold with $H$ being a degenerated distribution function at $\sigma^2$ and 
	\begin{align*}
	h_1(t)=\frac{\sigma^2}{\sigma^2+t},~h_2(t)=\frac{\sigma^4}{(\sigma^2+t)^2}.
	\end{align*}
	By the Stieltjes transformation of the MP law  \citep[see][Section 3.3]{bai2010spectral}, we have,
	\begin{align*}
	m_0(-\lambda)&=\frac{
		\sqrt{(\sigma^2-y \sigma^2+\lambda)^2+4y \lambda \sigma^2}-(\sigma^2-y \sigma^2+\lambda)}{2y \lambda \sigma^2},\\
	m'_0(-\lambda)&=\frac{\lambda (1+y)+\sigma^2(1-y)^2
	}{2y \lambda^2	\sqrt{(\sigma^2-y \sigma^2+\lambda)^2+4y \lambda \sigma^2}}-\frac{1-y}{2y \lambda^2}.
	\end{align*}
	Together with the fact that $R_1(\lambda)=\sigma^2 m_0(-\lambda)$ and $R_2(\lambda)=\sigma^4 m'_0(-\lambda)$, we would have
	\begin{align*}
	&\r^c_{\text{RLDA}}(\lambda) \inp\\
	&\Phi \left\{-\frac{ \Delta^2}{2 \sqrt{\Delta^2 +y_1+y_2}} \sqrt{\frac{2	\sqrt{(1+y+\lambda \sigma^{-2})^2-4y}}{\sqrt{(1+y +\lambda \sigma^{-2})^2-4y}+(1+y+\lambda \sigma^{-2})}}\right\}.
	\end{align*}	
We remark that Corollary 3.3 of \cite{dobriban2015high} presents the same result for a special case where $y=1$ and $\sigma^2=1$. Here in our results
we explicitly describe the effect of dimension and the tuning parameter $\lambda$ for general $y$ and $\sigma^2$. Figure \ref{fig1} show the simulation results for $y=0.5, y=1$ and $y=2$ from which we can see the empirical results and the theoretical conclusions are consistent.

\subsection{Isotropic case}
 Given the eigenvectors $v_1,\cdots,v_p$, we consider three different isotropic cases
 \begin{align}
 & \mu_1-\mu_2 \propto \frac{1}{\sqrt{p}}(v_1+\cdots+v_p), \label{iso1} \\
 & \Sigma^{-1/2}(\mu_1-\mu_2) \propto \frac{1}{\sqrt{p}}(v_1+\cdots+v_p), \label{iso2} \\
&  \Sigma^{-1}(\mu_1-\mu_2) \propto \frac{1}{\sqrt{p}}(v_1+\cdots+v_p).  \label{iso3}
 \end{align}
In statistics, these cases have different implications. \eqref{iso1} implies that the direct differences between the two population means, \eqref{iso2} indicates the differences after covariance standardization, and \eqref{iso3} implies that the true classification weights. Similar mean-covariance structures also arise in \cite{tony2014two} for the testing problem $H_0: \mu_1=\mu_2$.  In random matrix theory,  \cite{bai2007asymptotics} considered \eqref{iso1} when studying the asymptotic property of the eigenvectors and more details can be found in their Remark 1. 
 
To fix the Bayes error rate, we set $\Delta^2=(\mu_1-\mu_2)\trans \Sigma^{-1} (\mu_1-\mu_2)$ to be a constant and for each case, we list the results as follows.
\begin{itemize}
	\item Case \eqref{iso1}:
\begin{align*}
&\mu_1-\mu_2=\frac{\Delta}{\sqrt{\sum_{i=1}^p \lambda_i^{-1}}}(v_1+\cdots+v_p),\\
& h_1(t)=\left\{ \int \frac{1}{x}dH(x)\right\}^{-1} \int \frac{1}{x+t} dH(x), \\
& h_2(t)=\left\{ \int \frac{1}{x}dH(x)\right\}^{-1} \int \frac{x}{(x+t)^2} dH(x),\\
& H_1(\lambda)=\left\{ \int \frac{1}{x}dH(x)\right\}^{-1}m_0(-\lambda),\\
& H_2(\lambda)=\left\{ \int \frac{1}{x}dH(x)\right\}^{-1}\frac{m_0(-\lambda)-\lambda
	m'_0(-\lambda)}{\{1-y[1-\lambda m_0(-\lambda)]\}^2}.
\end{align*}	
	\item Case \eqref{iso2}:
\begin{align*}
&\mu_1-\mu_2=\frac{\Delta}{\sqrt{p}}(\sqrt{\lambda_1}v_1+\cdots+\sqrt{\lambda_p}v_p),\\
& h_1(t)=\int \frac{x}{x+t} dH(x),~h_2(t)=\int \frac{x^2}{(x+t)^2} dH(x),\\
& H_1(\lambda)=R_1(\lambda),~H_2(\lambda)=R_2(\lambda).
\end{align*}
	\item Case \eqref{iso3}:
\begin{align*}
&\mu_1-\mu_2=\frac{\Delta}{\sqrt{\sum_{i=1}^p \lambda_i}}(\lambda_1 v_1+\cdots+\lambda_p v_p),\\
& h_1(t)=\left\{ \int x dH(x)\right\}^{-1} \int \frac{x^2}{x+t} dH(x),\\
& h_2(t)=\left\{ \int x dH(x)\right\}^{-1} \int \frac{x^3}{(x+t)^2} dH(x),\\
& H_1(\lambda)= \frac{1}{1-y[1-\lambda m_0(-\lambda)]}-\left\{ \int x dH(x)\right\}^{-1} \frac{\lambda-\lambda^2 m_0(-\lambda)}{\{1-y[1-\lambda m_0(-\lambda)]\}^2},\\
& H_2(\lambda)=(1+R_1(\lambda))^2\left\{1-\frac{\lambda R_1(\lambda)}{\int x dH(x)} \right\}+\left\{ y-\frac{\lambda+2y \lambda R_1(\lambda)}{\int x dH(x)}\right\} R_2(\lambda).
\end{align*}
\end{itemize}

From the Bayesian perspective, we can use random effects that assume the true parameters are random to describe isotropic cases. For example, \cite{dobriban2015high} considered random weights:
\begin{align*}
\mu_1-\mu_2\ind (Z_1,\cdots,Z_p) \trans,
\end{align*}
where $Z_1,\cdots,Z_p$ are i.i.d random variables. Under their setting, 
	\begin{align*}
& (\mu_1-\mu_2)\trans \Sigma^{-1/2}(I_p+t\Sigma^{-1})^{-1}  \Sigma^{-1/2} (\mu_1-\mu_2)  \as c_1 \int \frac{1}{x+t} dH(x),\\
& (\mu_1-\mu_2)\trans \Sigma^{-1/2}(I_p+t\Sigma^{-1})^{-2}  \Sigma^{-1/2} (\mu_1-\mu_2) \as c_2\int \frac{x}{(x+t)^2} dH(x),
\end{align*}
and then $H_k(\lambda)$ have similar forms as our case \eqref{iso1}. Similarly, we can also consider the random effects
\begin{align*}
~\Sigma^{-1/2}(\mu_1-\mu_2) \ind (Z_1,\cdots,Z_p) \trans,~\mbox{or}~\Sigma^{-1}(\mu_1-\mu_2)  \ind (Z_1,\cdots,Z_p) \trans,
\end{align*}
and the results would exactly correspond to the cases \eqref{iso2} and \eqref{iso3}, respectively.

\subsection{Sparse case}
\begin{figure}[!htbp]
	\centerline{
		\begin{tabular}{ccc}				
			\psfig{figure=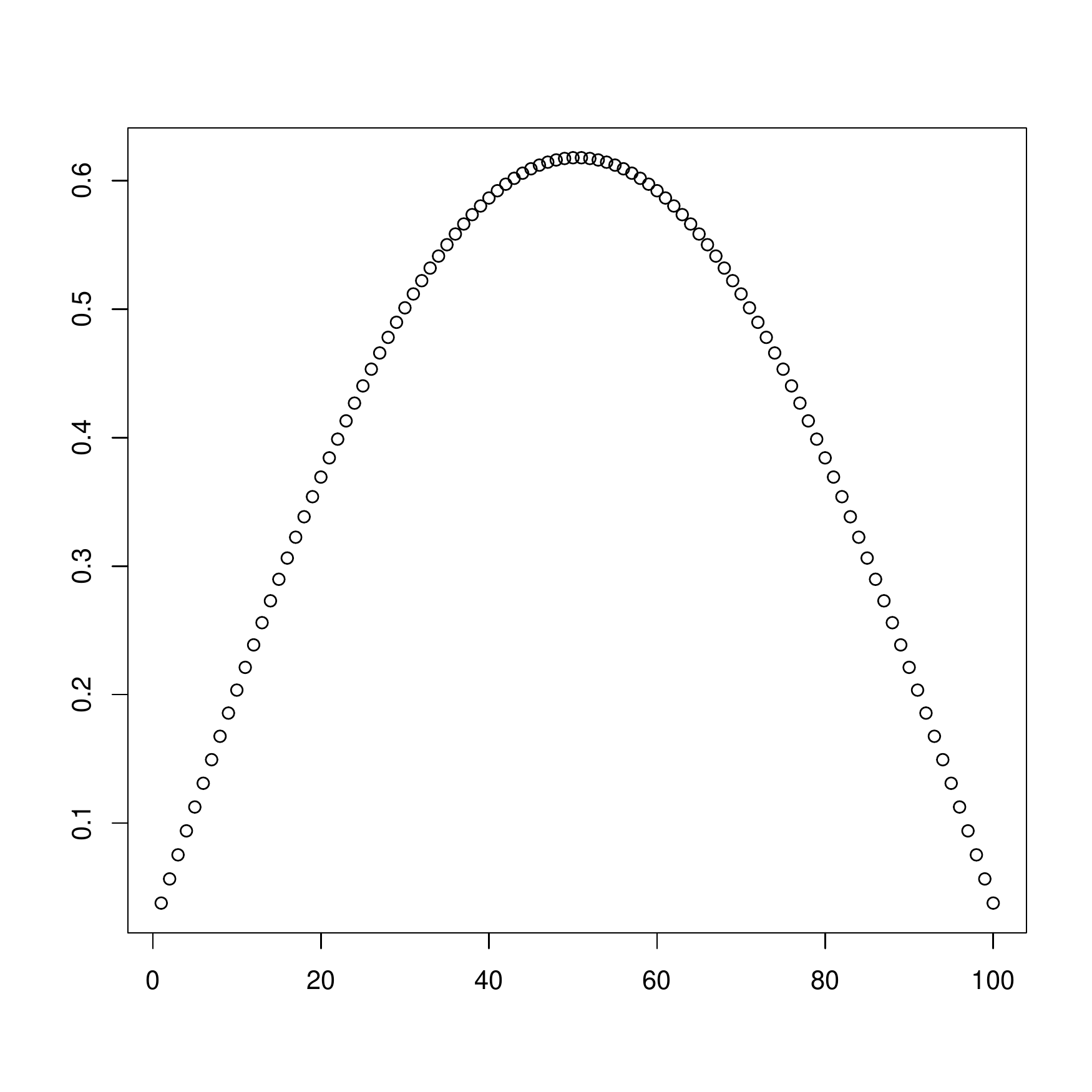,width=1.6 in,height=1.7 in,angle=0} &
			\psfig{figure=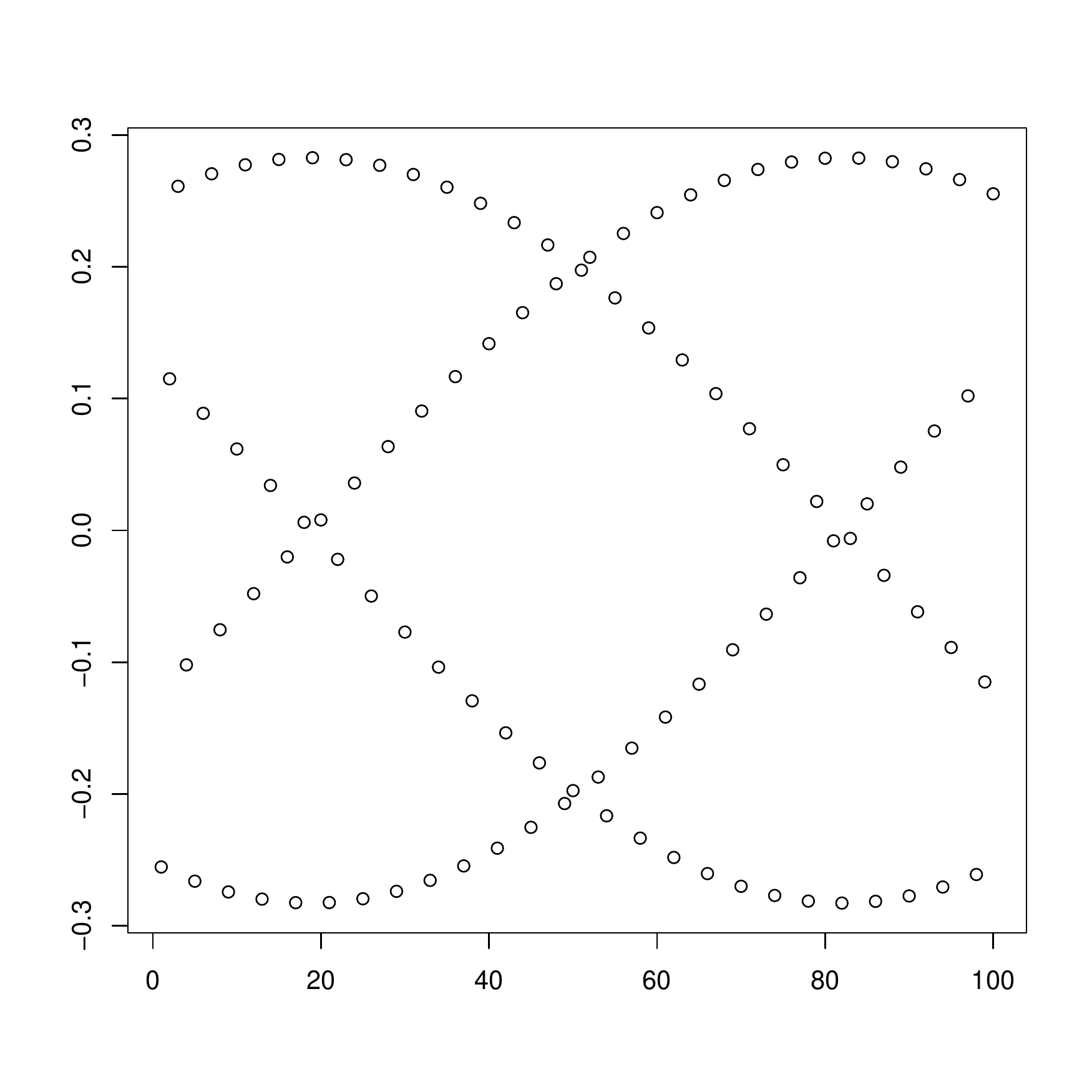,width=1.6 in,height=1.7 in,angle=0} &
			\psfig{figure=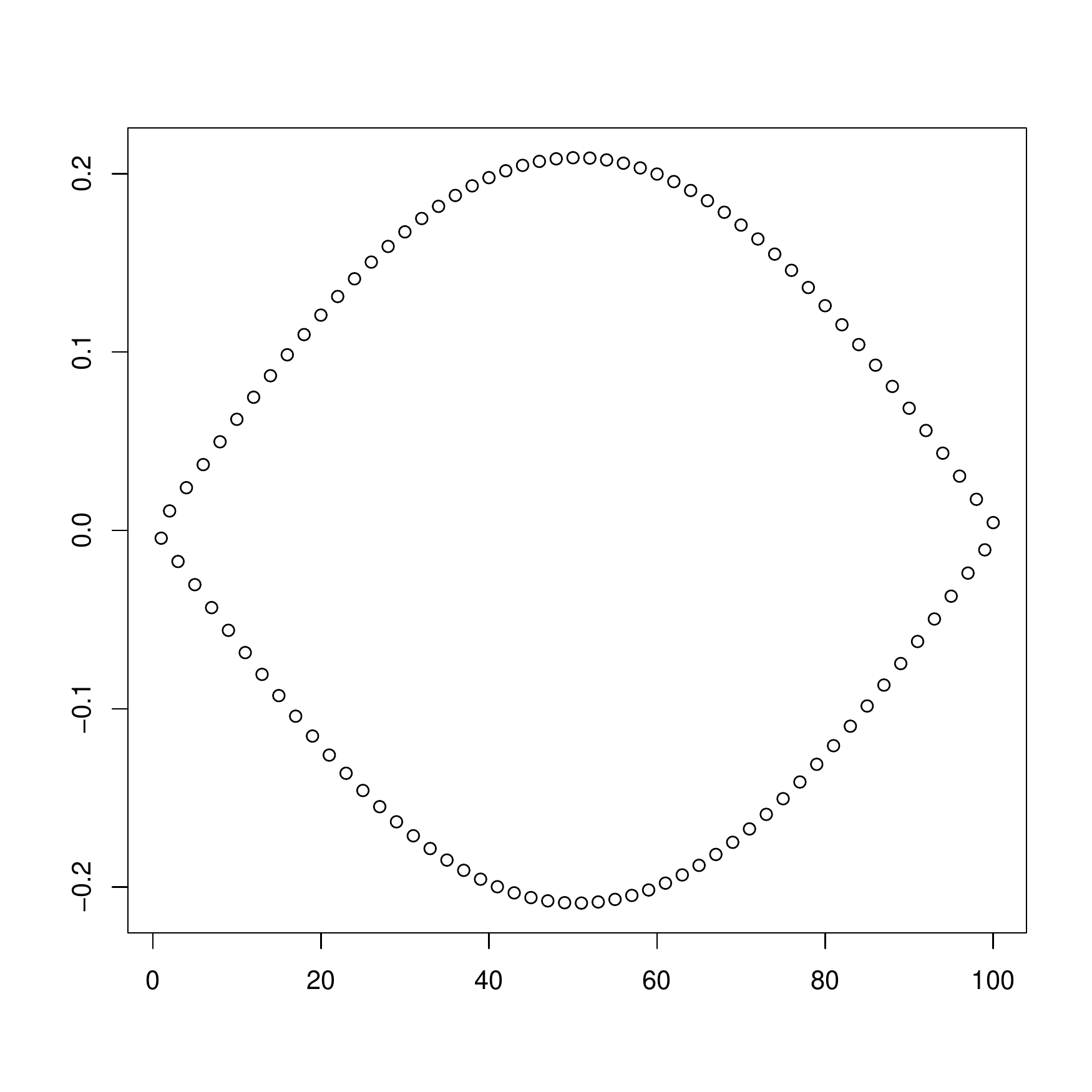,width=1.6 in,height=1.7 in,angle=0}\\
			$v_1$ & $v_{50}$ & $v_{100}$
		\end{tabular}
	}
	\caption{The eigenvectors of $\Sigma=(0.5^{|i-j|})_{100 \times 100}$ corresponding to $\lambda_1$,~$\lambda_{50}$ and $\lambda_{100}$.}
	\label{fig2}
\end{figure}

\begin{figure}[!htbp]
	\centerline{
		\begin{tabular}{c}				
			\psfig{figure=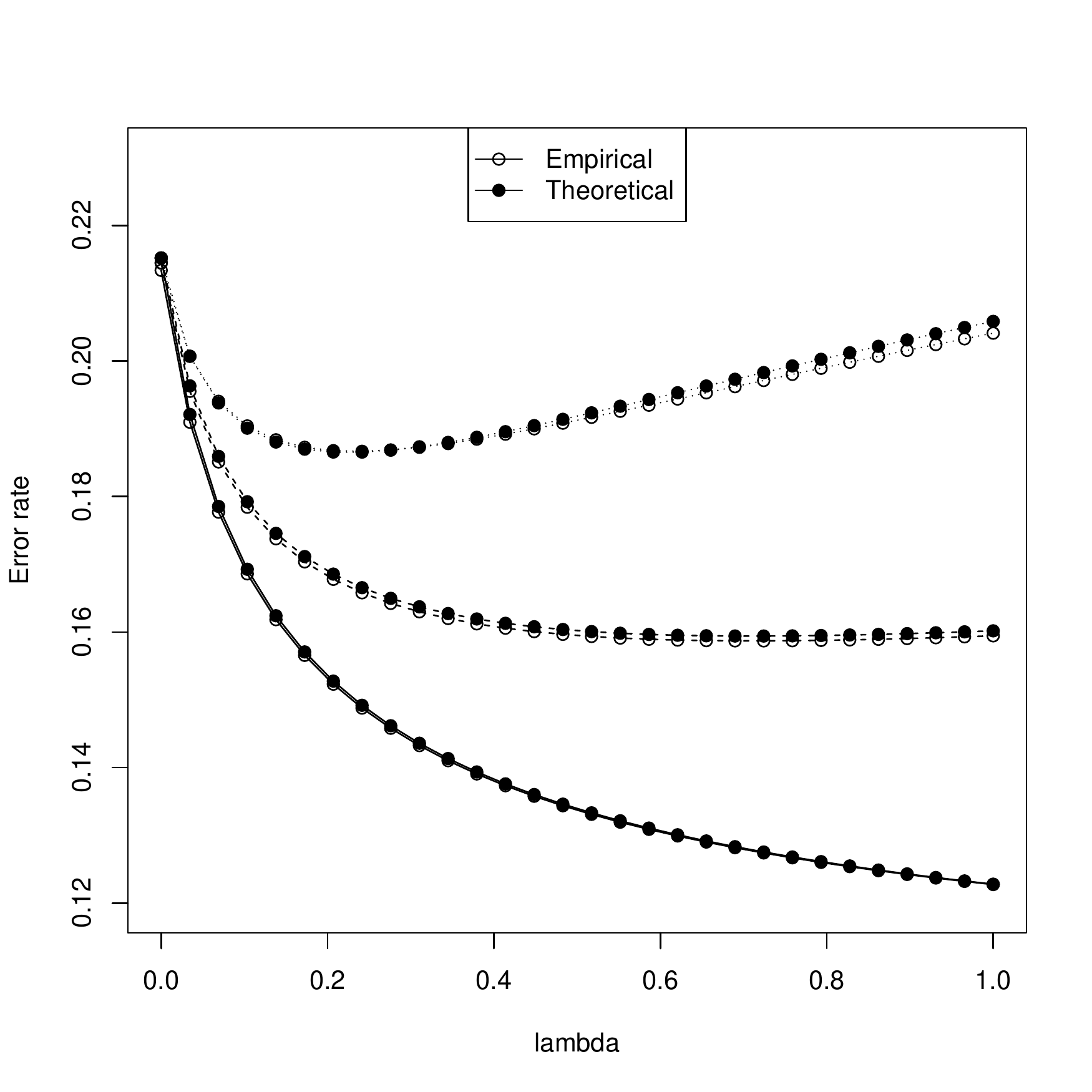,width=3.0 in,height=3.0 in,angle=0}
		\end{tabular}
	}
	\caption{The empirical and theoretical misclassification rates for  $\Sigma=(0.5^{|i-j|})_{100 \times 100}$ and $n_1=n_2=100$. The line stands for $\mu_1-\mu_2=4.435\cdot v_1$; the dashed line is the results for $\mu_1-\mu_2=2.005 \cdot v_{50}$ and the dotted line is the one for $\mu_1-\mu_2=1.480 \cdot v_{100}$. For all the cases, the Bayes error rate is $10 \%$.}
	\label{fig3}
\end{figure}

In this subsection, we consider the case where
\begin{align*} 
\mu_1-\mu_2 \propto v_k.
\end{align*}
We view this case as a sparse case since $\mu_1-\mu_2 \propto (0,\cdots,0,1,0,\cdots,0)$ under the basis vectors $\{v_1,\cdots,v_p\}$. When $\mu_1-\mu_2 \propto v_k$, we have $\Sigma^{-1}(\mu_1-\mu_2) \propto v_k$, indicating that the optimal linear classification direction is parallel to the eigenvector $v_k$.  We consider $\mu_1-\mu_2=\Delta \sqrt{\lambda_i} v_i$. Then
\begin{align*}
& \Delta^{-2}(\mu_1-\mu_2)\trans \Sigma^{-1/2}(I_p+t\Sigma^{-1})^{-1}  \Sigma^{-1/2} (\mu_1-\mu_2)=\frac{\lambda_i}{\lambda_i+t},  \\
& \Delta^{-2} (\mu_1-\mu_2)\trans \Sigma^{-1/2}(I_p+t\Sigma^{-1})^{-2}  \Sigma^{-1/2} (\mu_1-\mu_2)=\frac{\lambda_i^2}{(\lambda_i+t)^2}.
\end{align*}
When $\lambda_i$ has a limit which we still denote as $\lambda_i$, we have
\begin{align*}
h_1(t)=\frac{\lambda_i}{\lambda_i+t}, h_2(t)=\frac{\lambda_i^2}{(\lambda_i+t)^2},
\end{align*}
and 
\begin{align*}
  H_1(\lambda)=&\frac{1}{1-y+y\lambda m_0(-\lambda)+\lambda/\lambda_i},\\
  H_2(\lambda)=&\frac{1}{\{1-y+y\lambda m_0(-\lambda)+\lambda/\lambda_i\}^2}\\
 &\left\{\frac{1}{1-y(1-\lambda
	m_0(-\lambda))}-\frac{y(\lambda m_0(-\lambda)-\lambda^2
	m'_0(-\lambda))}{(1-y(1-\lambda m_0(-\lambda)))^2}\right\}.
\end{align*}
To illustrate the results, we consider a toy example $\Sigma=(\rho^{|i-j|})_{p \times p}$ with $|\rho|<1$. The covariance matrix $\Sigma$ is related to a stationary AR(1) process and is also used in \cite{bickel2004} for LDA. By the Szeg\"o theorem, we have, 
\begin{align*}
\lambda_k \approx \frac{1-\rho^2}{1+\rho^2-2 \rho \cos{\frac{k \pi}{p+1}}}.
\end{align*} 
Thus, $\lambda_1 \to (1+\rho)/(1-\rho),~\lambda_p \to (1-\rho)/(1+\rho)$ and $\lambda_{[p/2]} \to (1-\rho^2)/(1+\rho^2)$. Figure \ref{fig2} show the corresponding eigenvectors and Figure \ref{fig3} presents the results for $\mu_1-\mu_2 \propto v_k,~k=1,50,100$. We can see from Figure \ref{fig3} that the dimension effects have different performances for each $v_k$ and the tuning parameter $\lambda$ brings in different effects on the misclassification error rate. For all the simulations, the theoretical misclassification rates are consistent with the empirical ones.

\section{Simulations}
In this section, we conduct several simulations to illustrate the results. We consider the covariance matrix structures as follows:
\begin{align*}
\Sigma=(\rho^{|i-j|})_{p \times p},
\end{align*}
where $\rho$ reflects the correlations between the covariates and $|\rho|<1$. For all of our simulations, we let $\mu_2=0$ and rescale $\mu_1$ to control the Bayes error rate as $10 \%$. Here we fix the sample size $n=200$ and all the results are based on 100 replications. 
\subsection{Bias correction for LDA and RLDA}
\begin{figure}[!htbp]
	\centerline{
		\begin{tabular}{ccc}				
			\psfig{figure=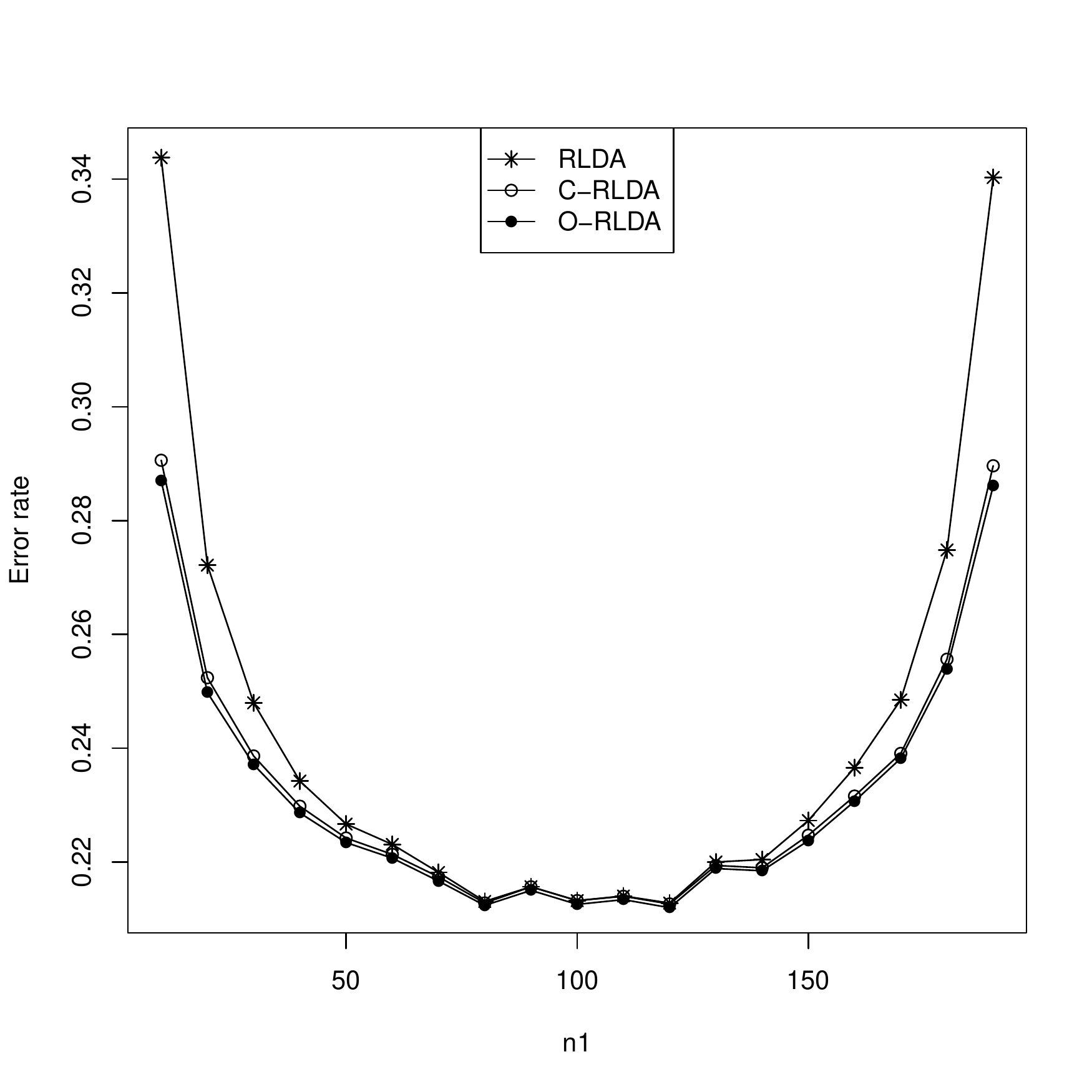,width=1.6 in,height=1.7 in,angle=0} &
		\psfig{figure=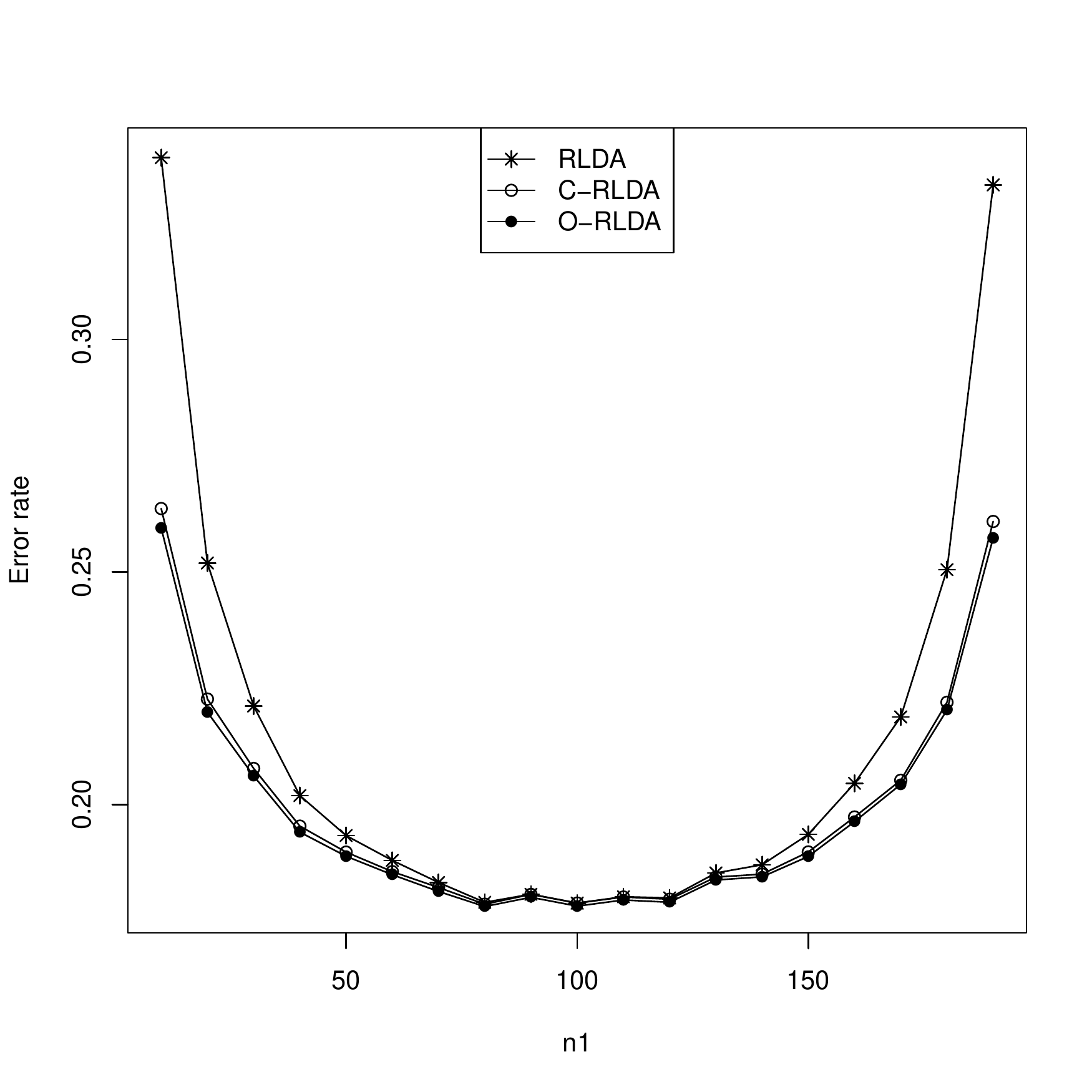,width=1.6 in,height=1.7 in,angle=0} &
		\psfig{figure=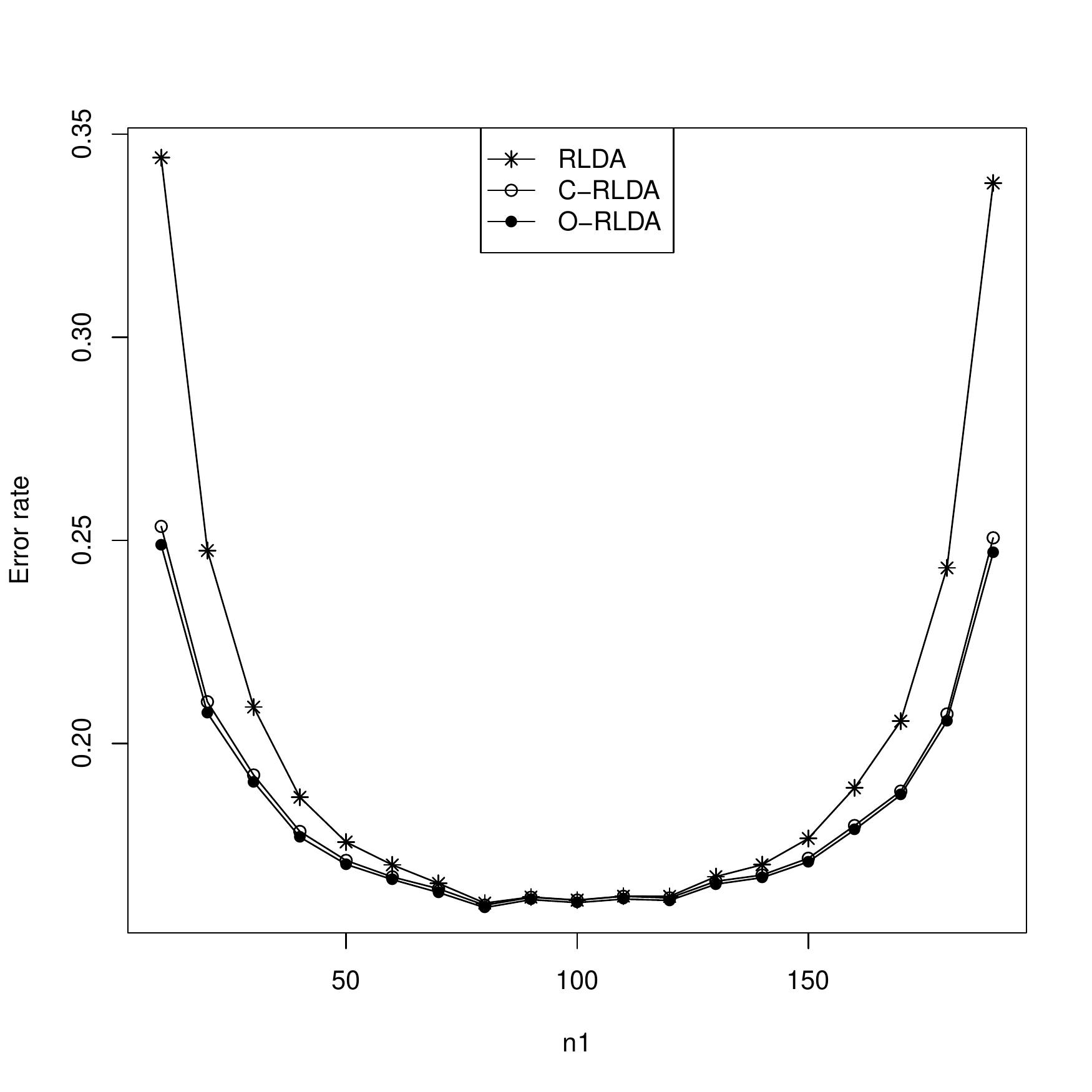,width=1.6 in,height=1.7 in,angle=0} \\
			$p=100,~\lambda=0$ & $p=100,~\lambda=0.1$ & $p=100,~\lambda=0.5$\\
						\psfig{figure=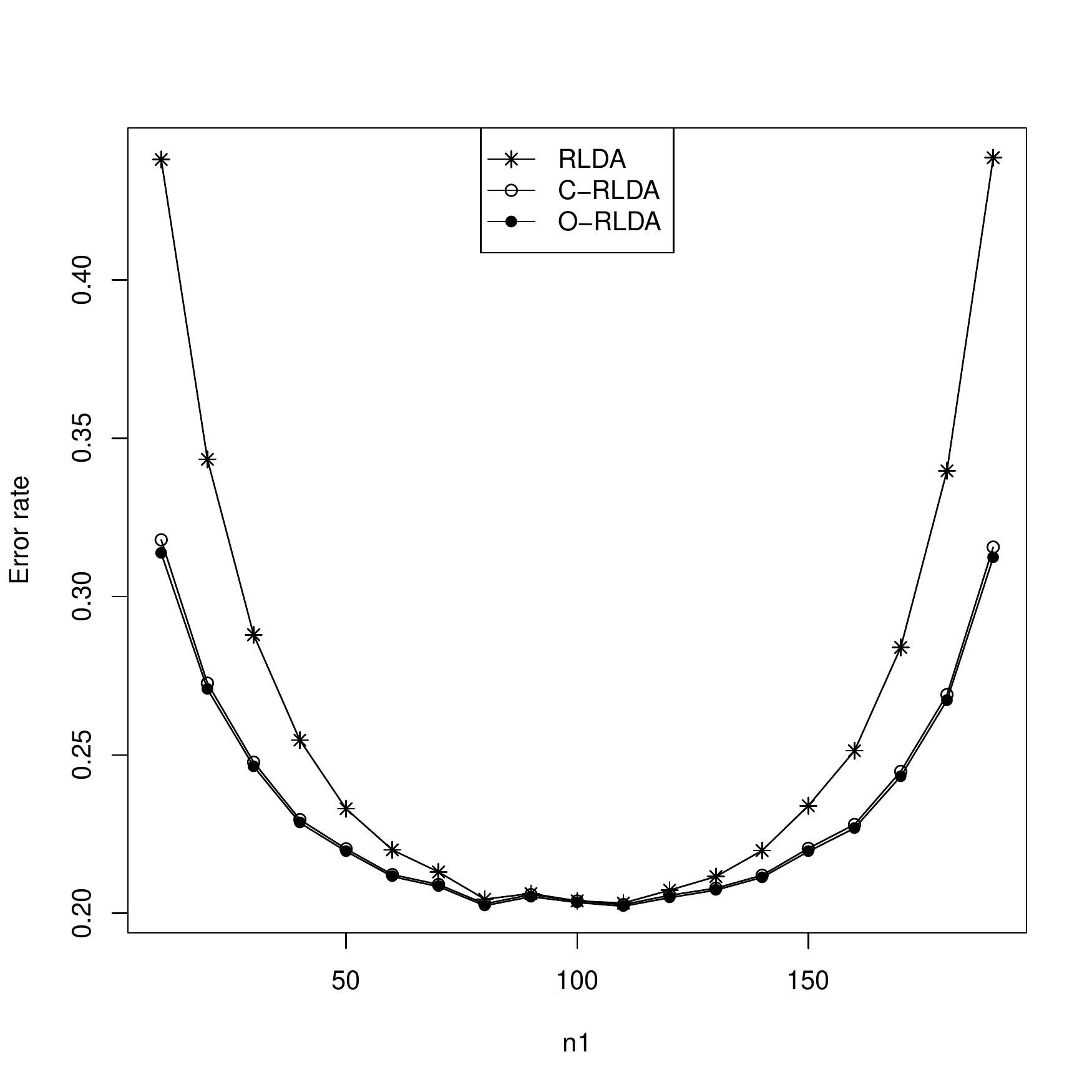,width=1.6 in,height=1.7 in,angle=0} &
			\psfig{figure=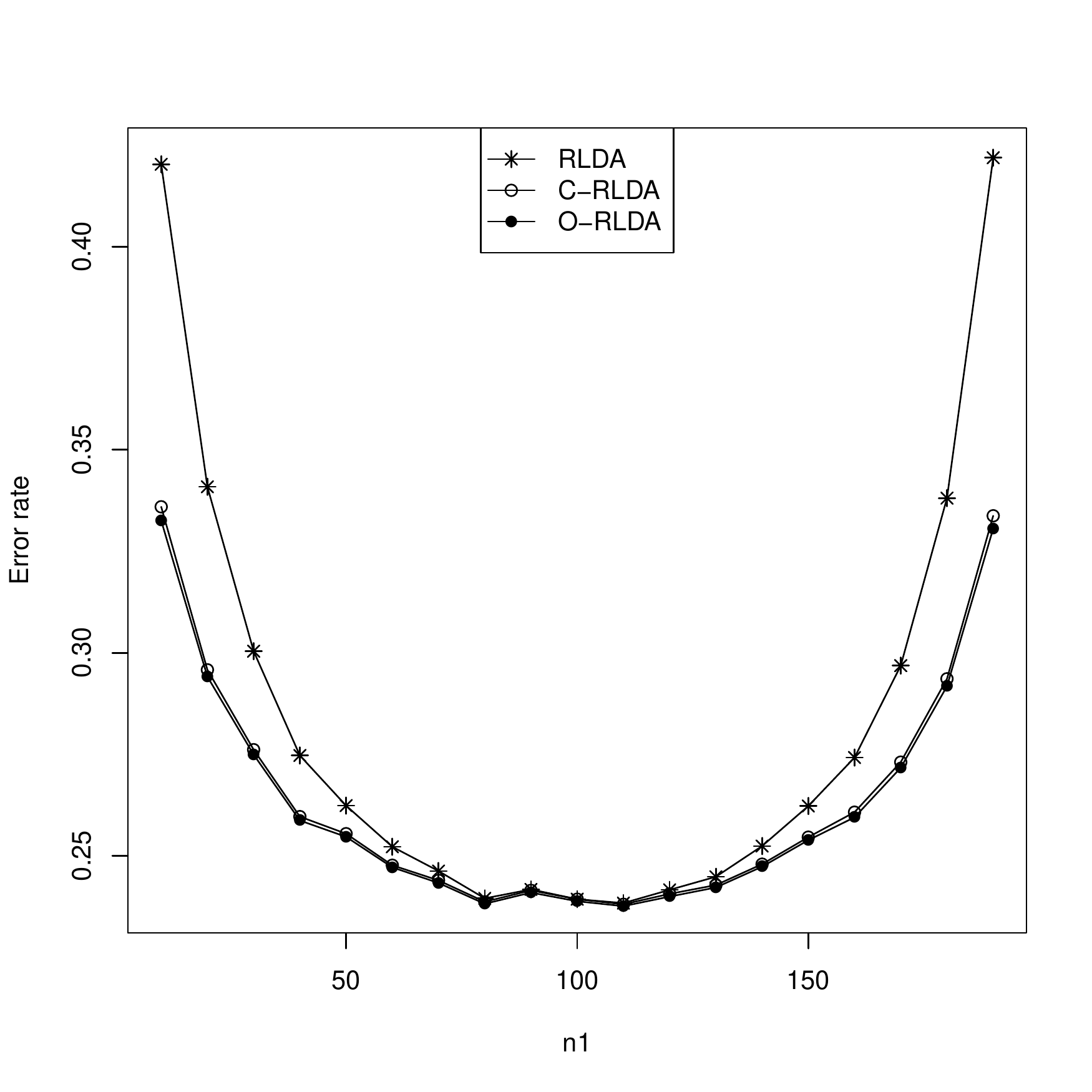,width=1.6 in,height=1.7 in,angle=0} &
			\psfig{figure=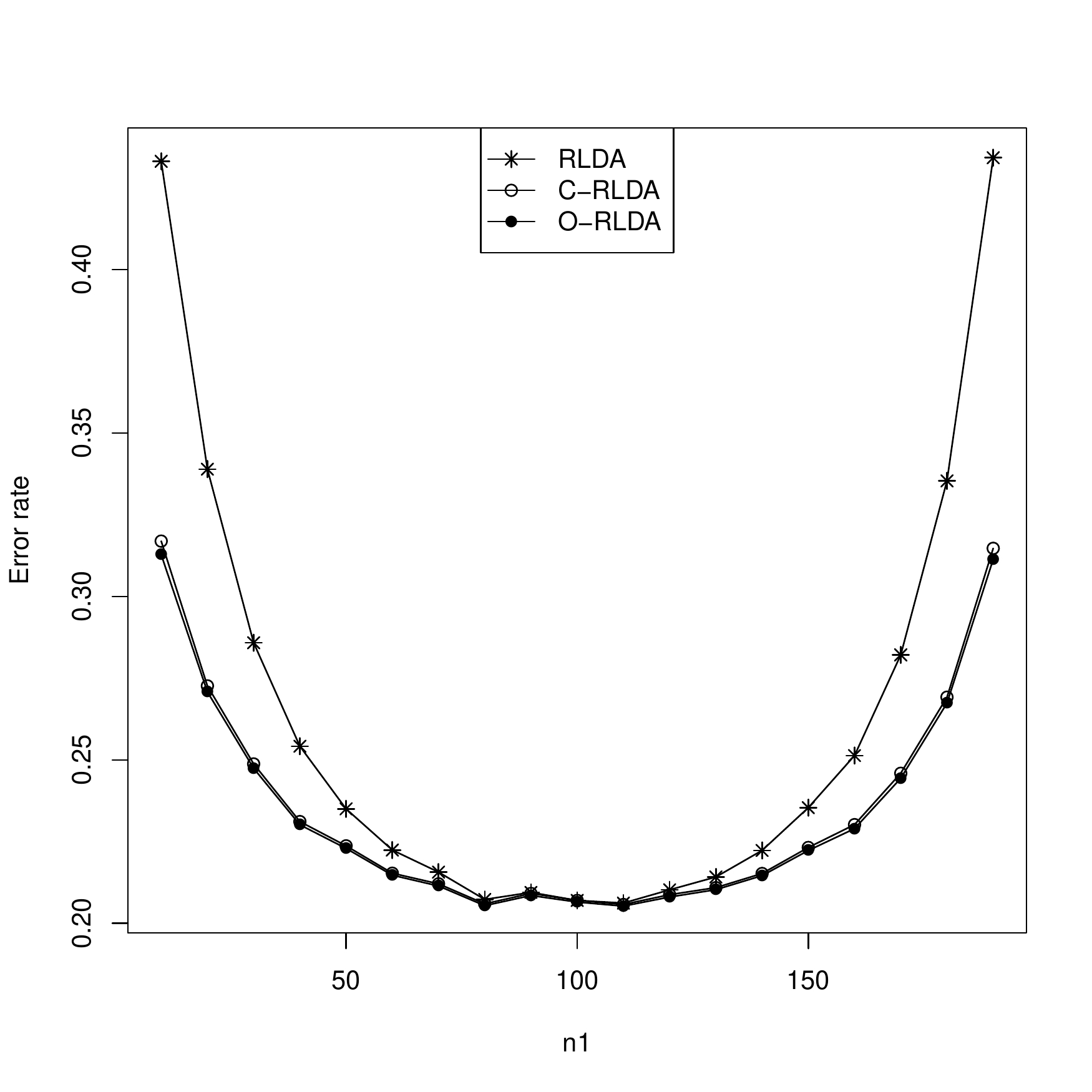,width=1.6 in,height=1.7 in,angle=0} \\
			$p=200,~\lambda=1$ & $p=200,~\lambda=0.1$ & $p=200,~\lambda=0.5$\\
						\psfig{figure=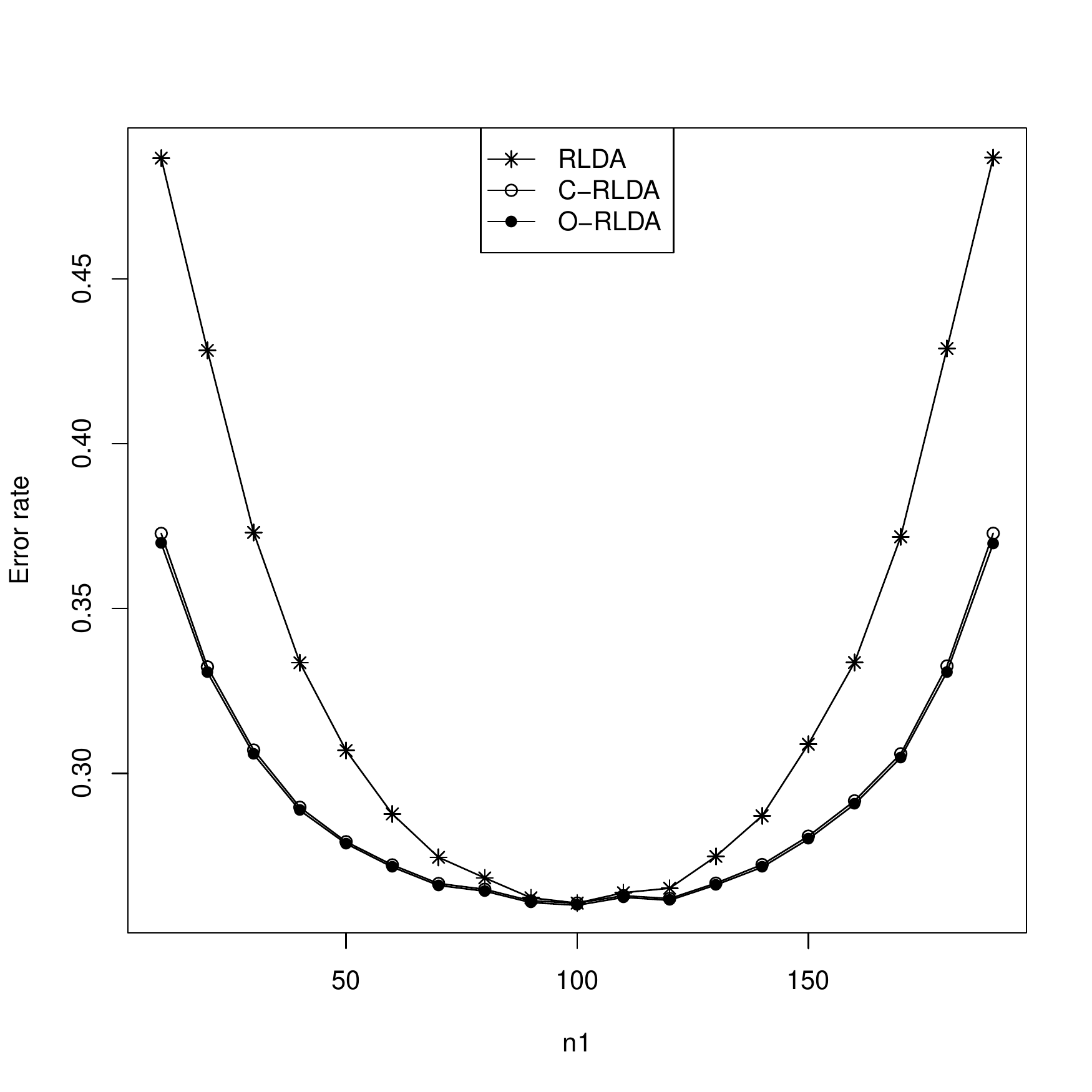,width=1.6 in,height=1.7 in,angle=0} &
			\psfig{figure=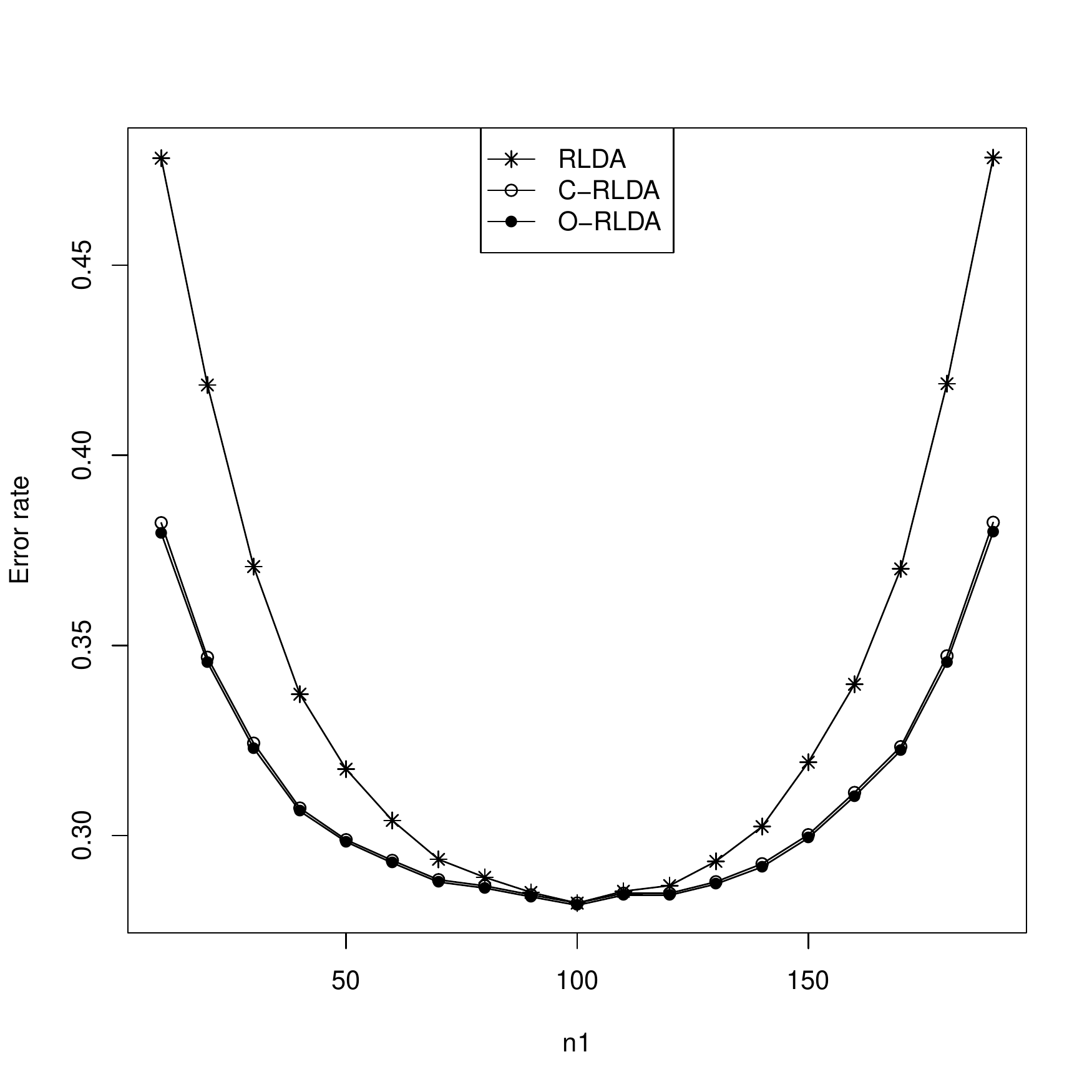,width=1.6 in,height=1.7 in,angle=0} &
			\psfig{figure=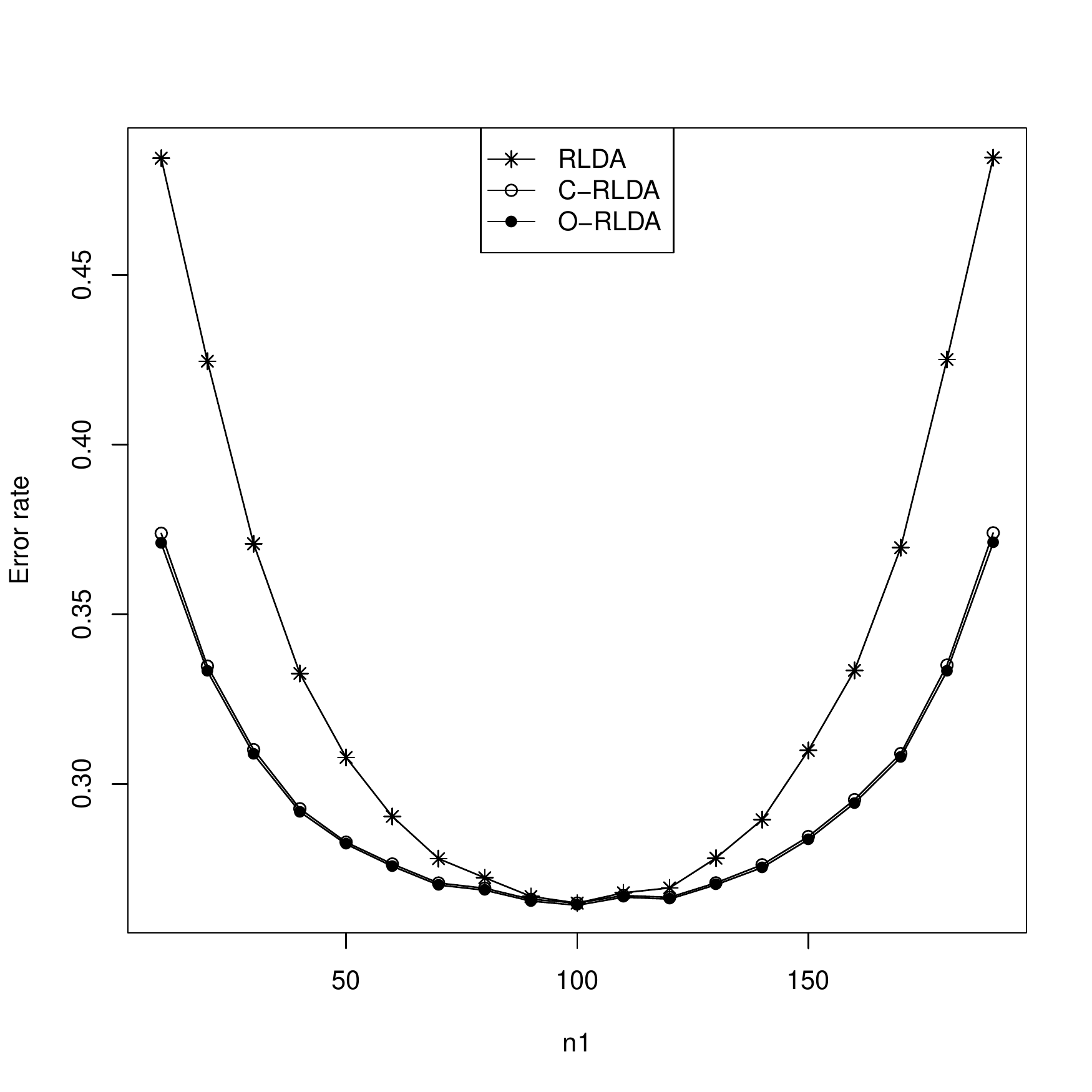,width=1.6 in,height=1.7 in,angle=0} \\
			$p=400,~\lambda=1$ & $p=400,~\lambda=0.1$ & $p=400,~\lambda=0.5$
		\end{tabular}
	}
	\caption{Simulations for RLDA, bias corrected RLDA (C-RLDA) and RLDA with optimal constant (O-RLDA).  The Bayes misclassification error rate is $10 \% $.}
	\label{fig4}
\end{figure}

When $n_1\neq n_2$, we propose a bias correction procedure to the intercept part for LDA and RLDA. Specially, we have
\begin{align*} 
\delta_{\text{LDA}}^c(X)=\mathbb{I}\left\{\Big(X-\frac{\bar{X}_1+\bar{X}_2}{2}\Big)\trans S_n^{-1}(\bar{X}_1-\bar{X}_2)+\frac{n-2}{n-p-3}\Big(\frac{p}{2 n_1}-\frac{p}{2 n_2}\Big)>0  \right\},
\end{align*}
and 
\begin{align*}
\delta^c_{\text{RLDA}}(X)=\mathbb{I}&\bigg\{\left(X-\frac{\bar{X}_1+\bar{X}_2}{2}\right)\trans (S_n+\lambda I_p)^{-1}(\bar{X}_1-\bar{X}_2) \nonumber \\
&+\left(\frac{p}{2n_1}-\frac{p}{2n_2}\right)\frac{1-\frac{1}{p}tr(\frac{1}{\lambda}S_n+I_p)^{-1}}{1-\frac{p}{n-2}+\frac{1}{n-2}tr(\frac{1}{\lambda}S_n+I_p)^{-1}}>0\bigg\},~\lambda>0.
\end{align*}
Our first simulation is to show the performance of the bias-corrected LDA and RLDA. We fix $n=200$ and let $n_1$ rang from 10 to 190. We set $\lambda$ as 0 (which corresponds to LDA),  0.1 and 0.5 for $p<n$ and 0.1, 0.5 and 1 when $p\geq n$ .The data dimension $p$ is set to be 100, 200 or 400. As a benchmark, we also include the linear classifier with optimal constant
 \begin{align*}
\alpha_r=-\frac{1}{2} (\mu_1+\mu_2) \trans(S_n+\lambda I_p)^{-1}(\bar{X}_1-\bar{X}_2),
\end{align*}
which is denoted by $\mbox{O-RLDA}$. We set $\Sigma=(0.5^{|i-j|})_{p \times p}$ and $ \mu_1 \propto (1,0,\cdots,0)$, and all the simulation results are plotted in Figure \ref{fig4}. We can observe that the bias corrected RLDA (C-RLDA) achieves less misclassification rate than the original RLDA and also the performance of C-RLDA is quite close the one of RLDA with optimal constants. Specially, the first figure of Figure \ref{fig4} is the results for LDA where $p=100,~n=200$. We also conduct simulations for more covariance and mean structures and the results follow similar patterns. 

\subsection{Dimension effect of LDA, naive Bayes and RLDA}
\begin{figure}[!htbp]
	\centerline{
		\begin{tabular}{ccc}				
			\psfig{figure=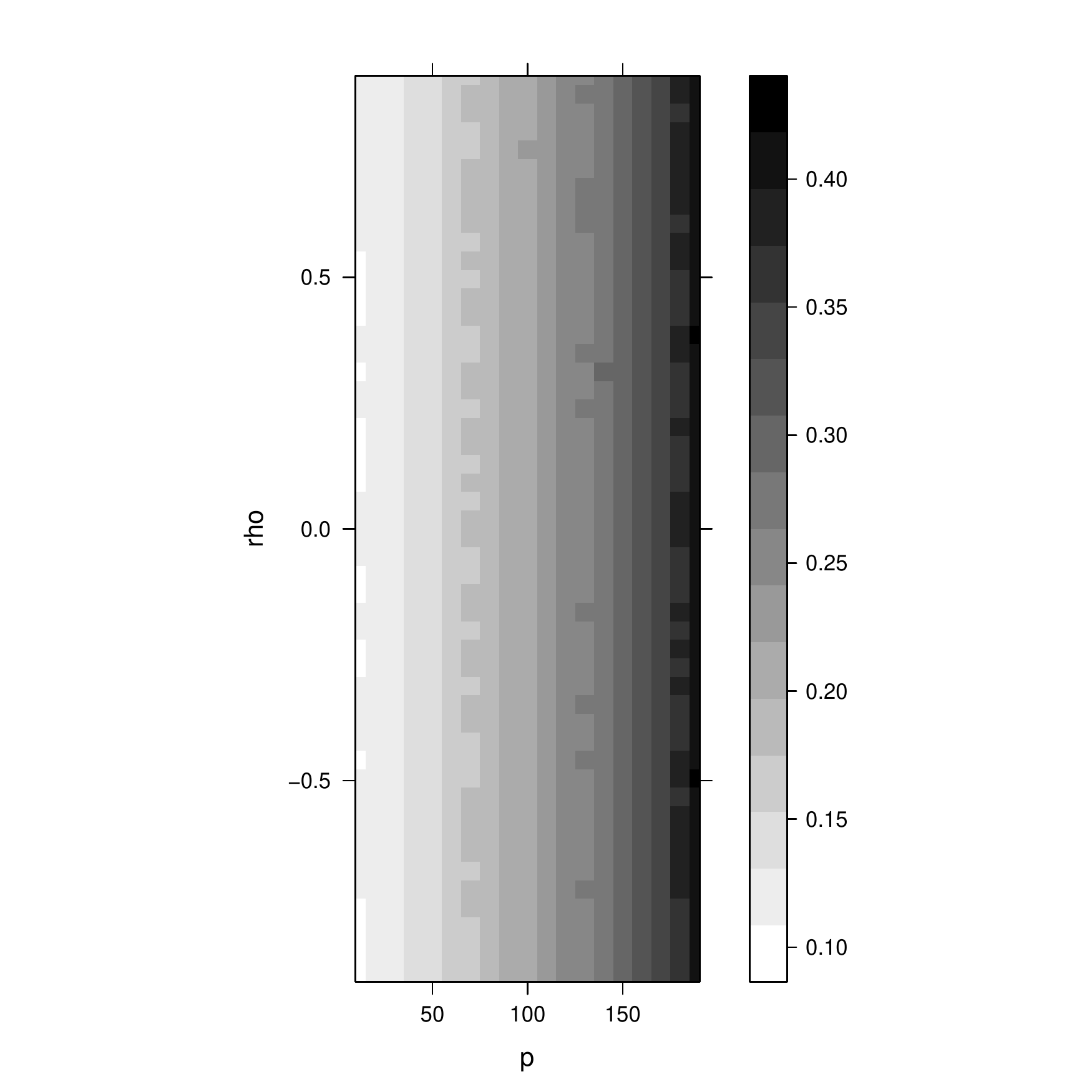,width=1.6 in,height=1.7 in,angle=0} &
			\psfig{figure=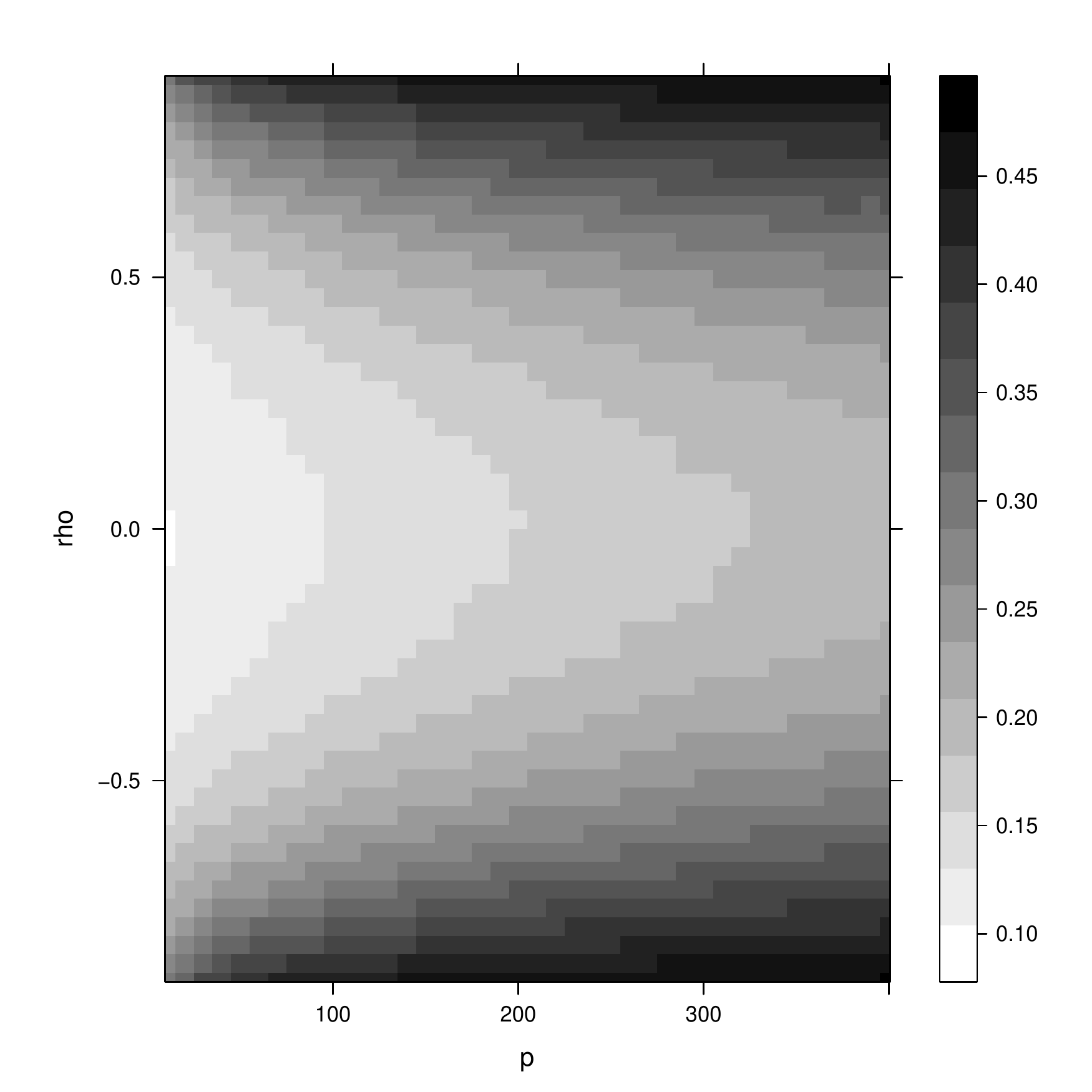,width=1.6 in,height=1.7 in,angle=0} &
				\psfig{figure=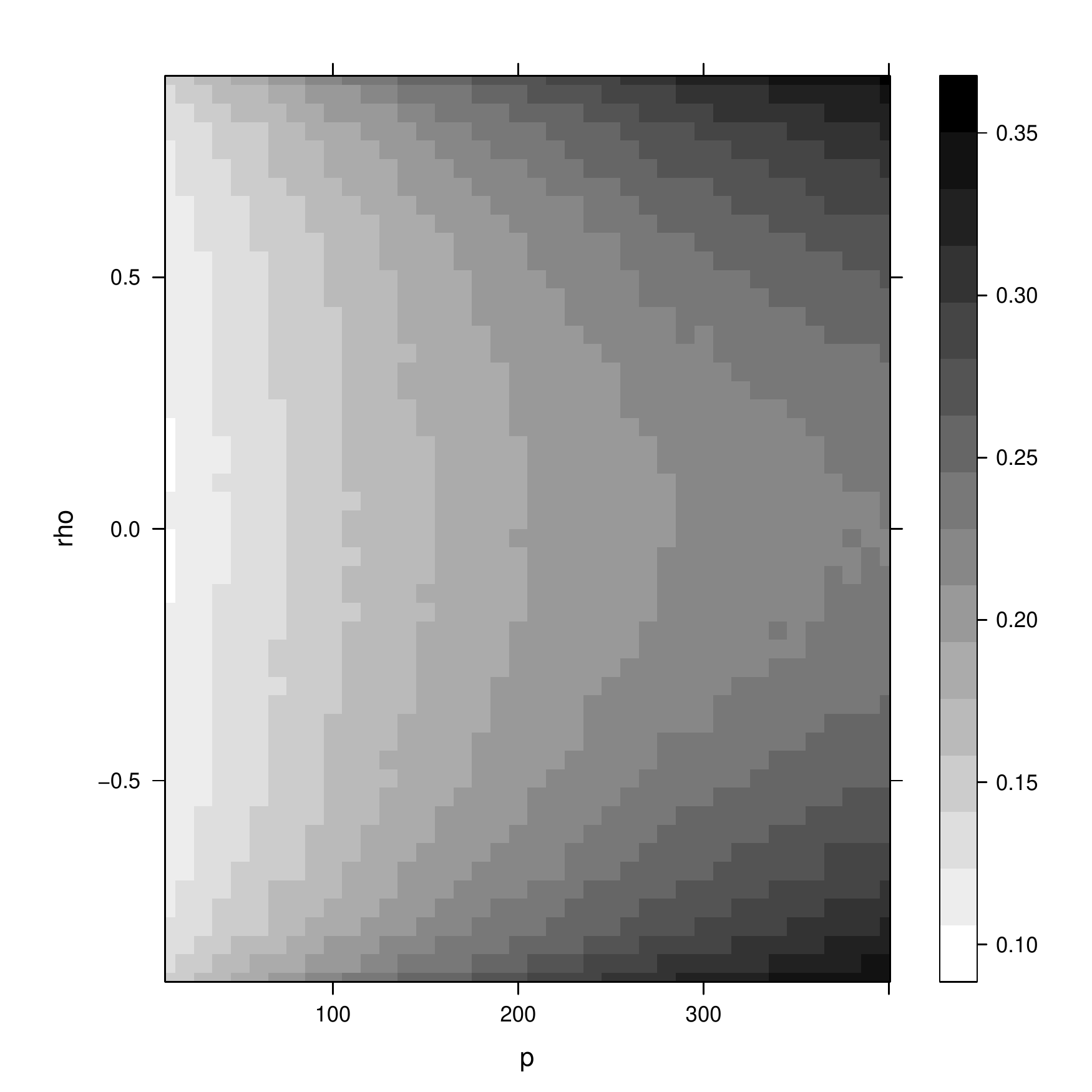,width=1.6 in,height=1.7 in,angle=0} \\
			LDA for Case 1  & Naive Bayes for Case 1 & RLDA with $\lambda=0.5$ for Case 1\\
					\psfig{figure=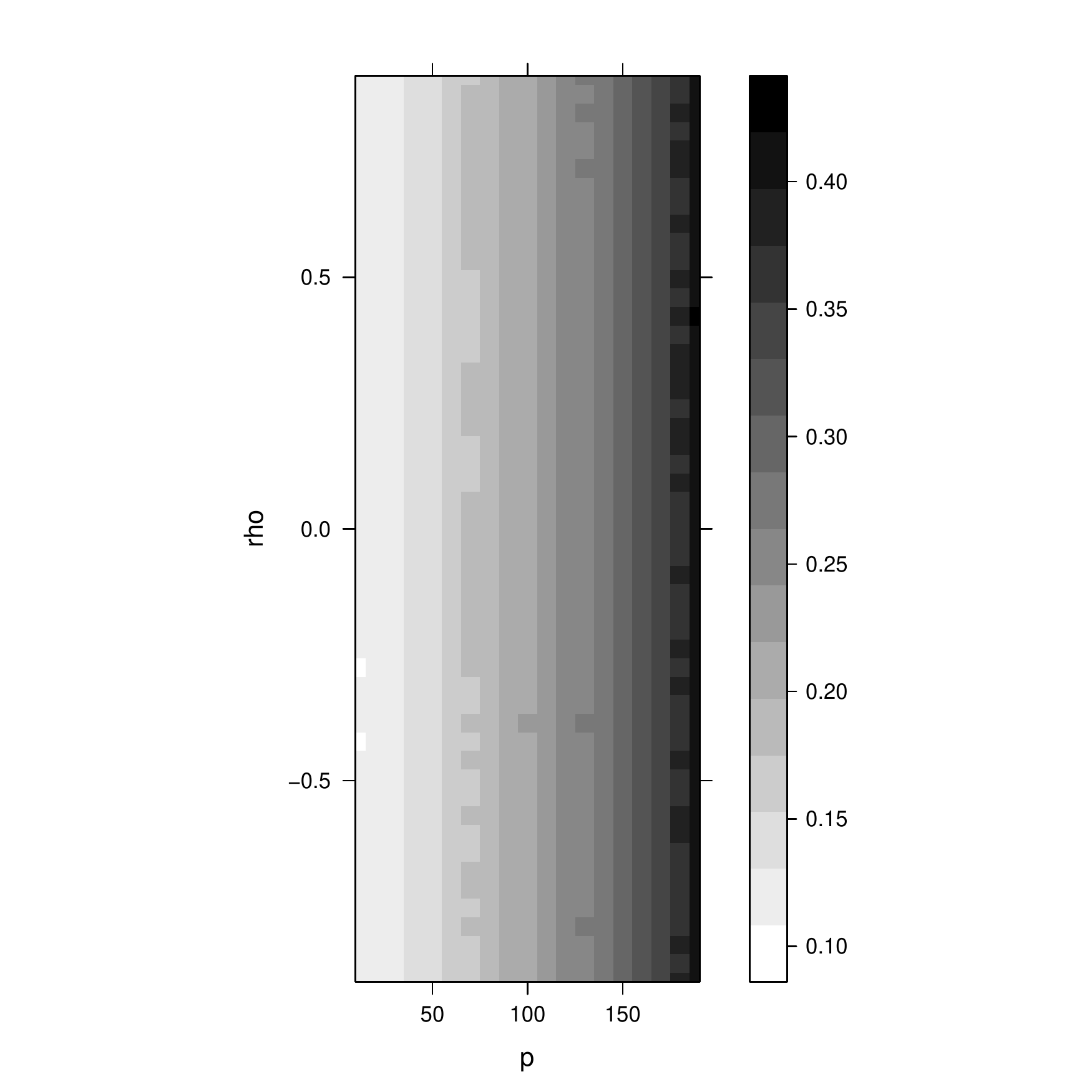,width=1.6 in,height=1.7 in,angle=0} &
		\psfig{figure=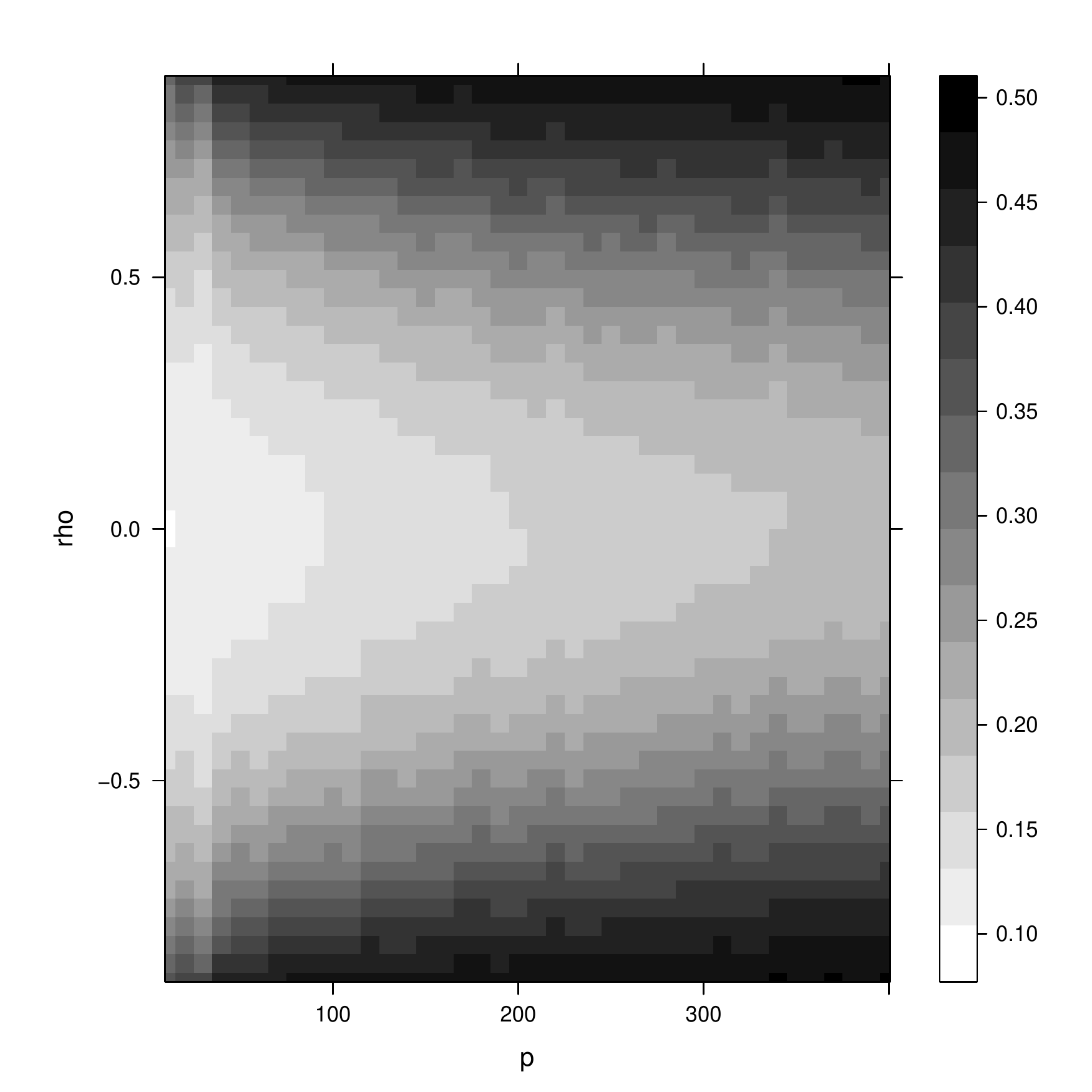,width=1.6 in,height=1.7 in,angle=0} &
		\psfig{figure=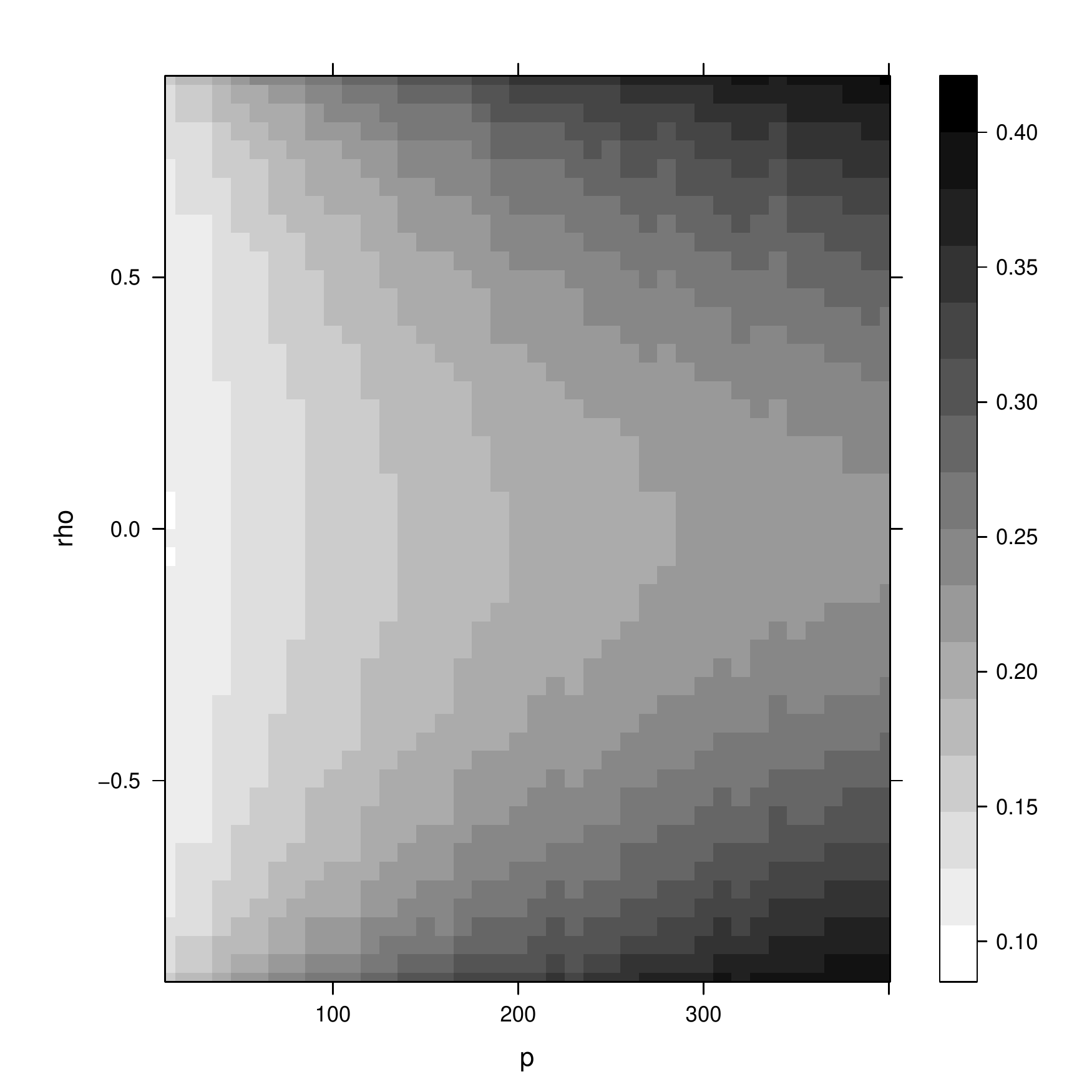,width=1.6 in,height=1.7 in,angle=0} \\
				LDA for Case 2& Naive Bayes for Case 2& RLDA with $\lambda=0.5$ for Case 2\\
					\psfig{figure=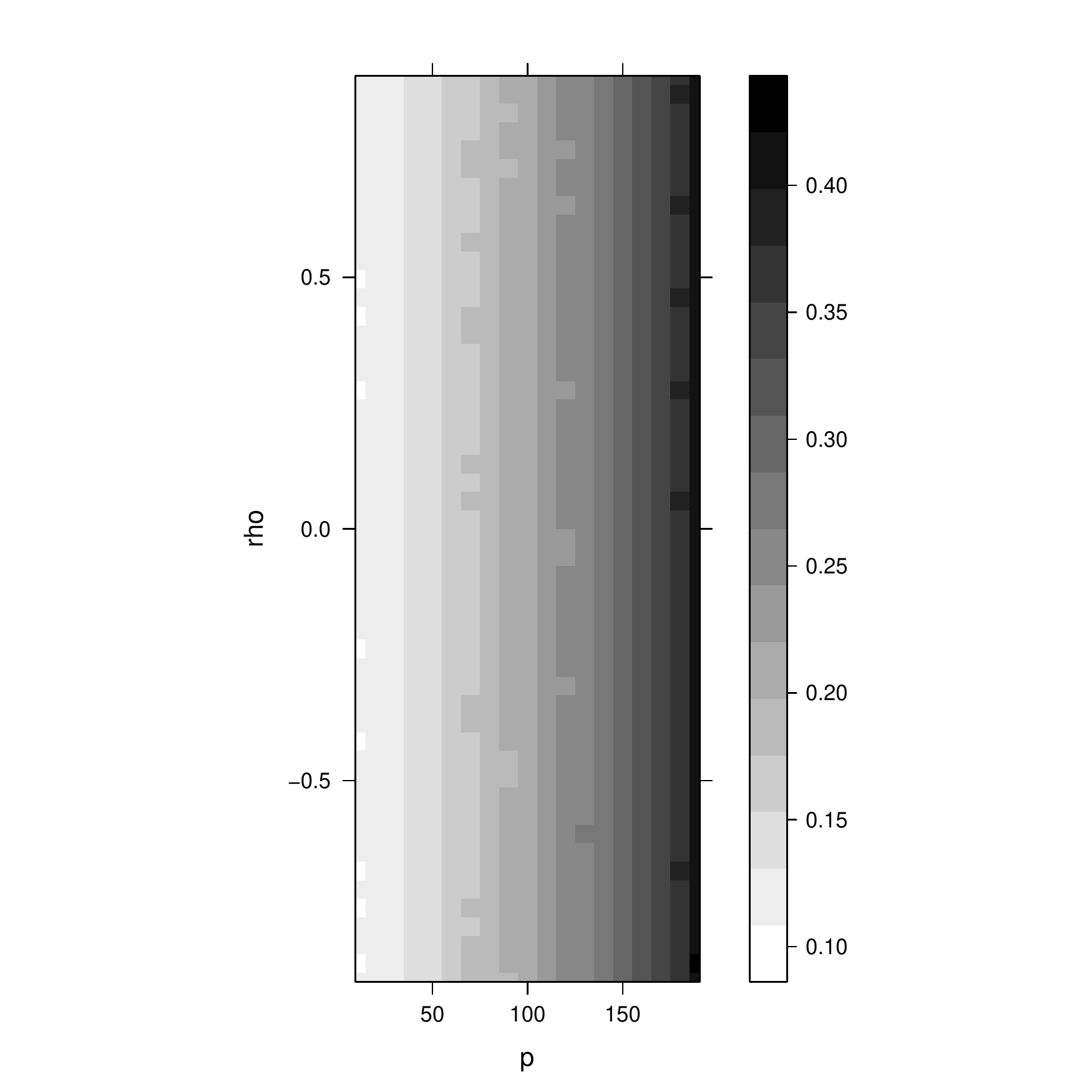,width=1.6 in,height=1.7 in,angle=0} &
		\psfig{figure=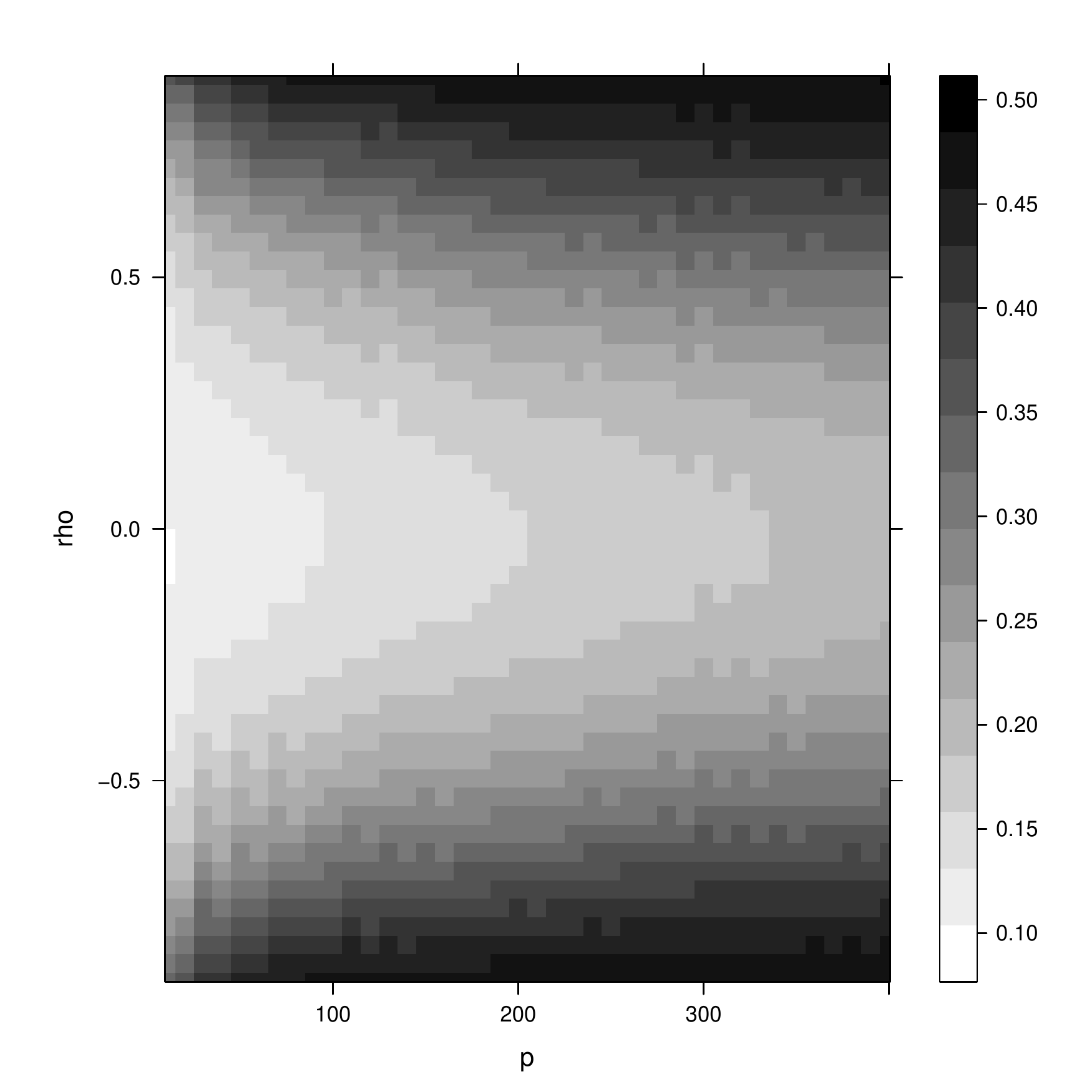,width=1.6 in,height=1.7 in,angle=0} &
		\psfig{figure=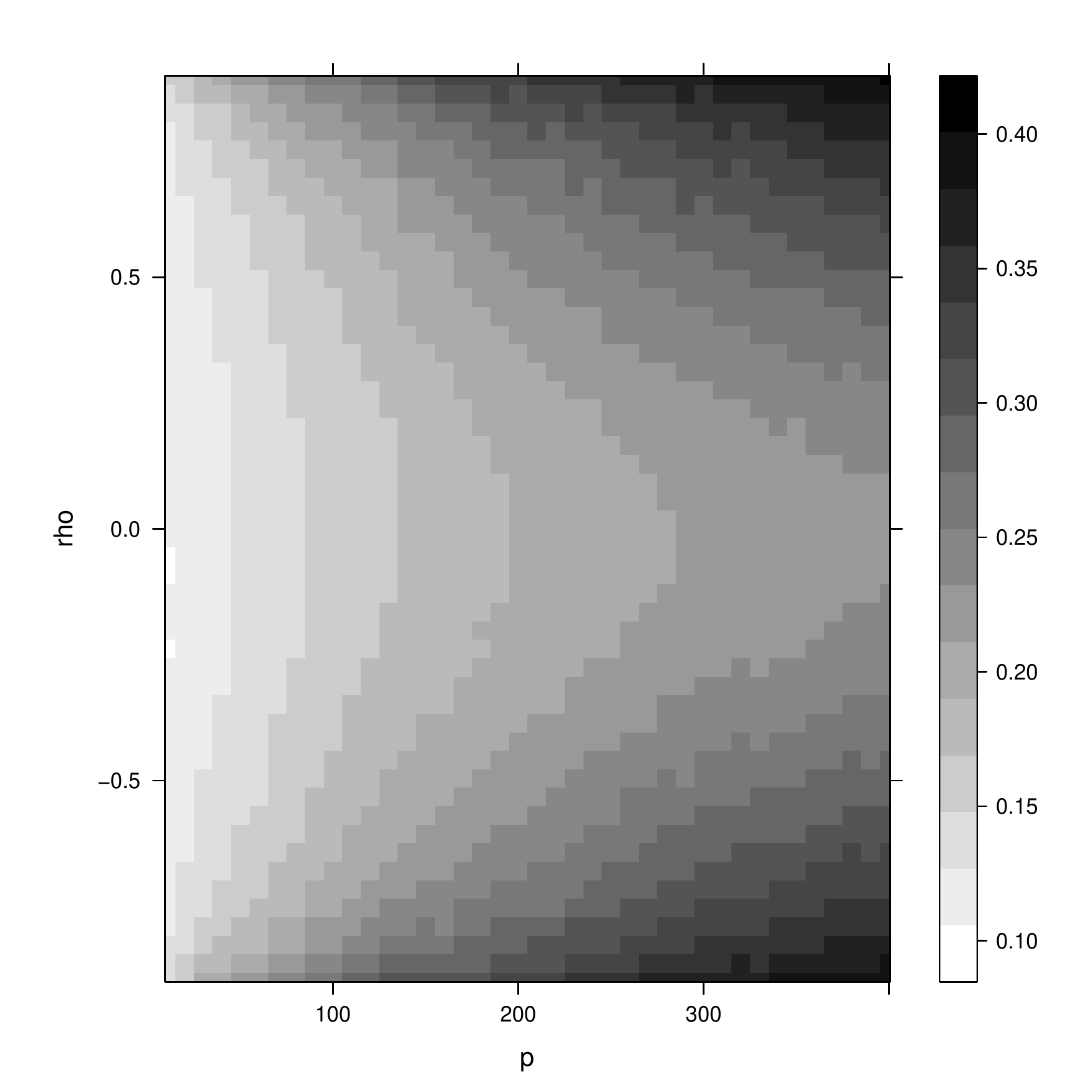,width=1.6 in,height=1.7 in,angle=0} \\
			LDA for Case 3& Naive Bayes for Case 3& RLDA with $\lambda=0.5$ for Case 3\\
		\end{tabular}
	}
	\caption{Heat maps of the error rates for LDA, Naive Bayes and RLDA with $\lambda=0.5$ where the horizontal is the data dimension $p$ ranging from 10 to 190 or 400 and the vertical is the correlation $\rho$ ranging from -0.9 to 0.9.   The Bayes misclassification error rate is $10 \% $.}
	\label{fig5}
\end{figure}

In this work, our main focus is to explicitly derive the dimension effects for LDA and RLDA. In detail, the asymptotic misclassification rate of LDA depends on the ratios $p/n_1,p/n_2$ and the one of RLDA also involves the structure of the covariance the means. More details can be found in Proposition \ref{prop0} and \ref{prop1}. In Remark \ref{nbayes}, we derive the dimension effect for naive Bayes. Here, we conduct a comprehensive comparison for LDA, naive Bayes and RLDA. For $\mu_1$, we consider three scenarios:
\begin{itemize}
	\item[] Case 1: $ \mu_1 \propto (1,0,\cdots,0)$;
	\item[] Case 2:  10 \%  elements of $\mu_1$ are from $N(0,1)$;
	\item[] Case 3: all the elements of $\mu_1$ are from $N(0,1)$;
\end{itemize}
 and for each case we still rescale $\mu_1$ to control the Bayes error rate to be 10\%.  All the results are presented in Figure \ref{fig5}. For each heat map of the error rate,  the horizontal is the data dimension $p$ ranging from 10 to 190 for LDA or 400 for naive Bayes and RLDA. The vertical is the correlation $\rho$ ranging from -0.9 to 0.9.  When the data dimension $p$ is increased, the misclassification rates for all the classifiers also increases which is well understood. When the correlations $|\rho|$ ranges from 0 to 0.9,  the performance of LDA almost has no changes and Naive Bayes gets worse and worse. For RLDA, the performance is affected by $|\rho|$, but is less sensitive than naive Bayes.

\subsection{RLDA and sparse LDA}
\begin{figure}[!htbp]
	\centerline{
		\begin{tabular}{ccc}				
		\psfig{figure=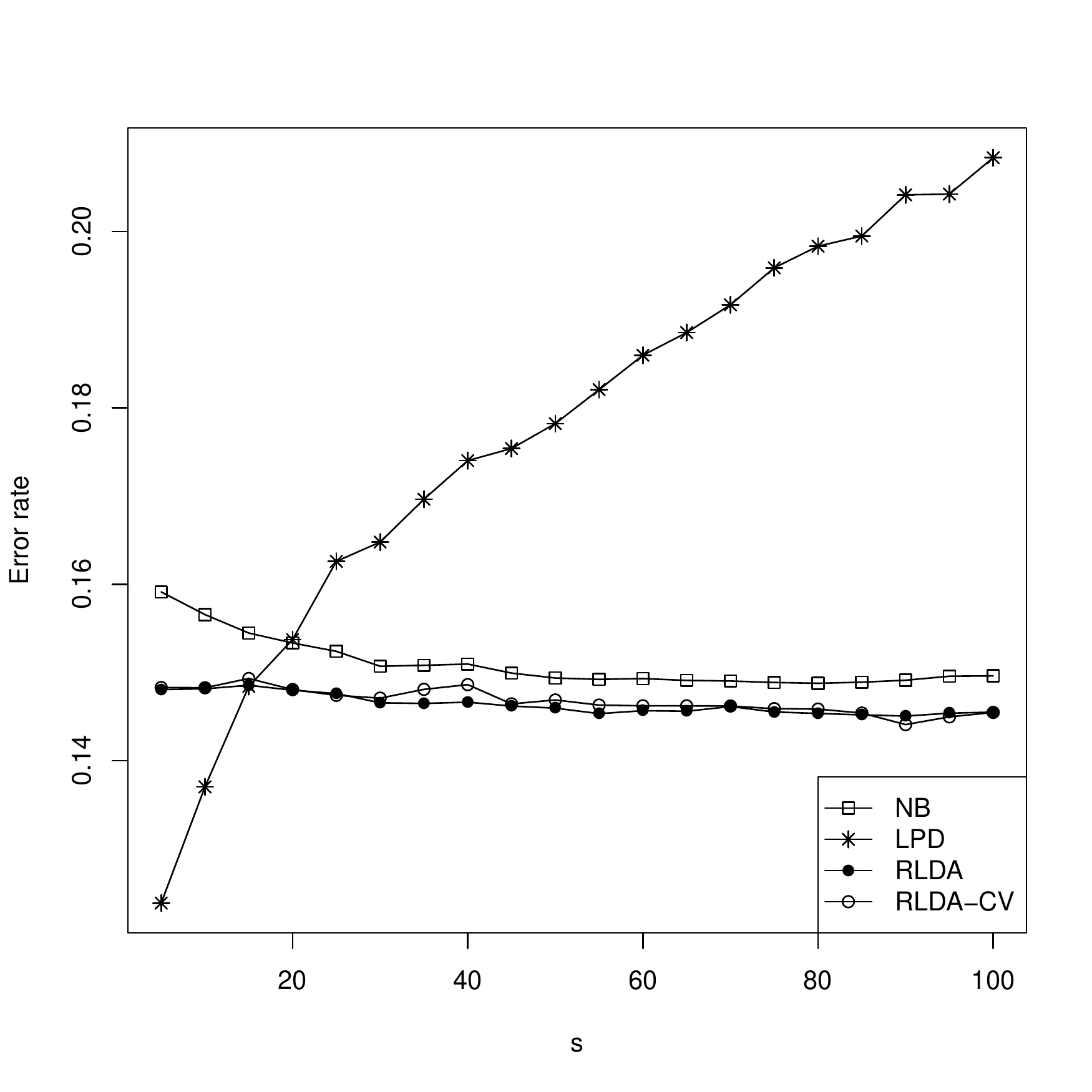,width=1.6 in,height=1.7 in,angle=0} &
\psfig{figure=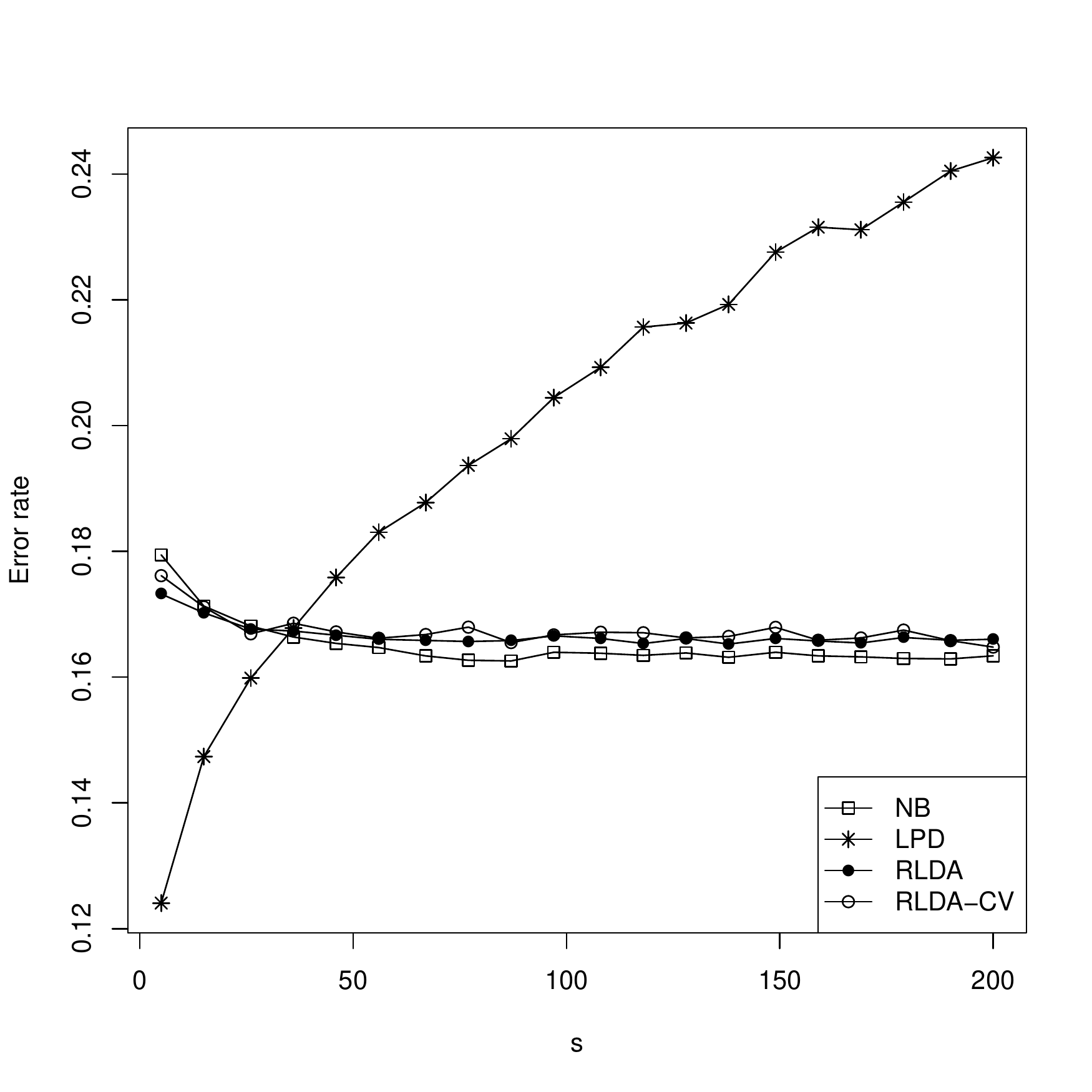,width=1.6 in,height=1.7 in,angle=0} &
\psfig{figure=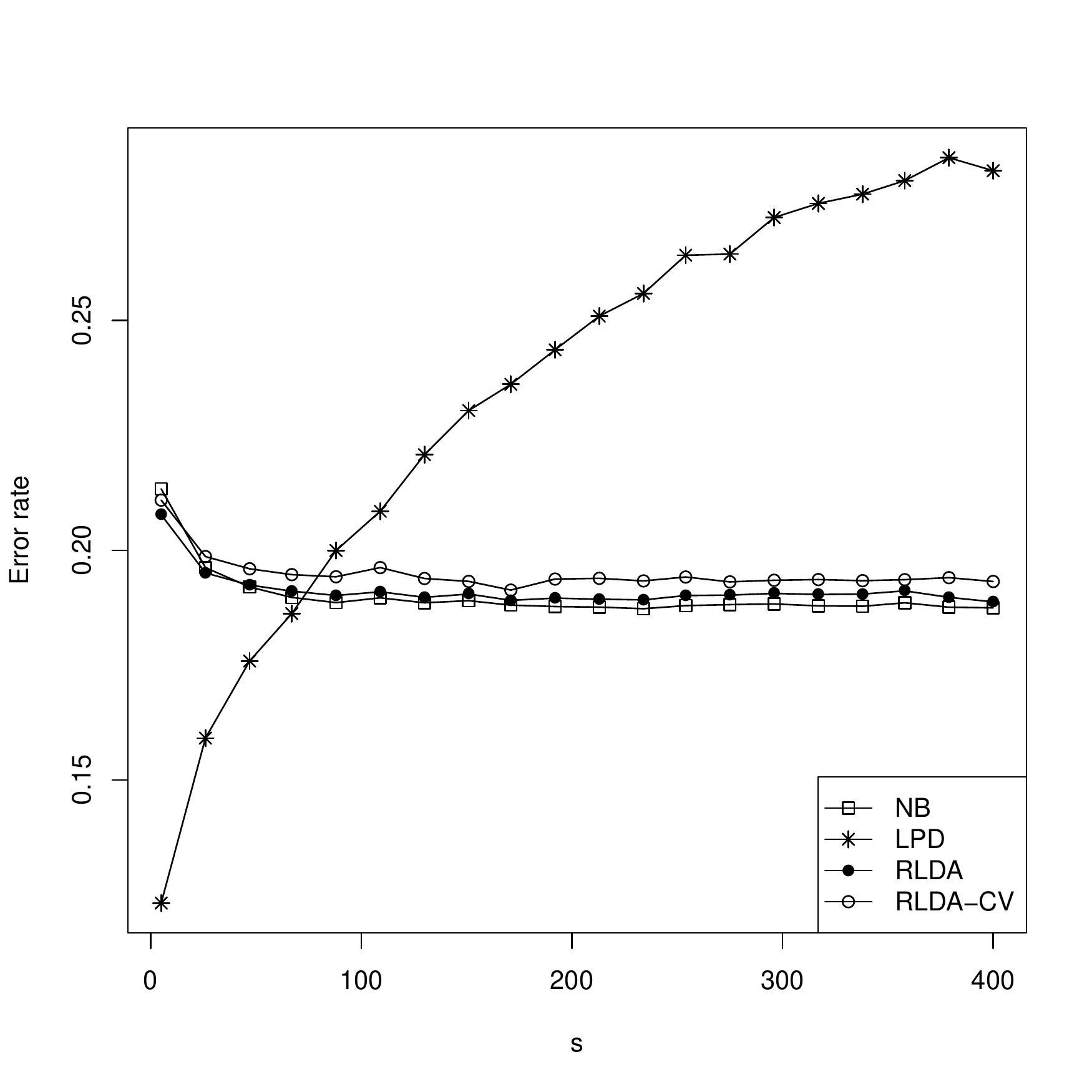,width=1.6 in,height=1.7 in,angle=0} \\
			$p=100$ & $p=200$& $p=400$\\
		\end{tabular}
	}
	\caption{Misclassification rates for naive Bayes, LPD and RLDA with cross validation and estimated $\lambda$. The horizontal is the number of non-zeros of $\Sigma^{-1}(\mu_1-\mu_2)$.   }
	\label{fig6}
\end{figure}

In this part, we compare RLDA with the LPD \citep{cai2011direct} which is one of sparse LDA methods. The tuning parameters are chosen by 5-folds cross validation.  We also include the RLDA where $\lambda_n$ is estimated by our method in Section \ref{slambda} and the naive Bayes. In summary, we have four methods NB, LPD, RLDA-CV and RLDA.  We control the sparsity level as
\begin{align*}
\Sigma^{-1}(\mu_1-\mu_2) \propto (\mbox{rnorm(s)},0,\cdots,0);
\end{align*}
with $\mbox{rnorm(s)}$ represent $s$ random variables generated independently from $N(0,1)$, and in this study $s$ is increased from 5 to $p$. Figure \ref{fig6} presents the simulation results. We observe that  when $s$ is small, LPD is outstanding with the least error rates and the performance becomes poorer and poorer  when $s$  increases, while RLDA and NB are robust to different $s$.  Furthermore, we find out that the version with estimated tuning parameter is comparable to the one with cross validation.  Overall, we claim that LPD is applicable to the sparse cases while RLDA is favorable to the dense cases. 

\section{Acknowledgments}
We thank two reviewers, an associate editor, and the editor for their most helpful comments. Wang is partially supported by Shanghai Sailing Program 16YF1405700 and National Natural Science Foundation of China 11701367. Jiang is partially supported by the Early Career Scheme from Hong Kong Research Grants Council PolyU 253023/16P.

\newpage
\section{Appendix}
\begin{lem} \label{lem1}
	Assuming $W_m \sim \mbox{Wishart}(I_p,m)$ where $m>p+7$, for any non-random unit  vector $e \in \mR^{p}$, we have
	\begin{align}
	& E(e \trans W_m^{-1} e)=\frac{1}{m-p-1}, \label{lem1a}\\
	& E(e \trans W_m^{-2} e)=\frac{m-1}{(m-p)(m-p-1)(m-p-3)},  \label{lem1b}
	\end{align}
	and 
	\begin{align}
	& E\{ (e \trans W_m^{-1} e)^2\}=\frac{1}{(m-p-1)(m-p-3)},  \label{lem1c}\\
	& E\{ (e \trans W_m^{-2} e)^2\}=\frac{m^2}{(m-p)^6} \Big\{1+O\Big(\frac{1}{m-p}\Big) \Big\} .  \label{lem1d}
	\end{align}	 
\end{lem}
\noindent\textit{Proof of Lemma \ref{lem1}:}  
For any non-random orthogonal matrix  $U$,  we have
\begin{align*}
U \trans   W_m U  \sim \mbox{Wishart}(I_p,m),
\end{align*}
and
\begin{align*}
e \trans W_m^{-1} e \ind (Ue) \trans W^{-1}_m U e,~e \trans W_m^{-2} e \ind (Ue) \trans W^{-2}_m U e.
\end{align*}
Setting $U e=(1,0,\cdots,0) \trans$ yields
\begin{align*}
& e \trans W_m^{-1} e \ind  (W_m^{-1})_{11},~e \trans W_m^{-2} e\ind (W_m^{-2})_{11}=\sum_{i=1}^p (W_m^{-1})^2_{1i}.
\end{align*}
Therefore, we can get the results \eqref{lem1a}-\eqref{lem1c} by \cite{von1988moments} or the Proposition 2.1 of \cite{cook2011mean}. The order of \eqref{lem1d} is derived based on Theorem 2 of  \cite{matsumoto2012general} and properties of the Weingarten function \citep{collins2006integration} which involves lengthy and tedious calculations and is hence neglected.   \hfill$\fbox{}$

\begin{lem}  \label{lem2}
	Assuming $Y \sim N(0,I_p)$ and $Y$ is independent with $S_n$, we have
	\begin{align*}
	 \frac{1}{p} Y \trans B_n(\lambda) Y \inp R_1(\lambda),~~{\rm and}~~ \frac{1}{p}Y \trans B^2_n(\lambda) Y \inp R_2(\lambda),
	\end{align*}
where $B_n(\lambda)$ and $R_1(\lambda), R_2(\lambda)$ are defined as in \eqref{Bn} and Theorem \ref{thm4} respectively. 
\end{lem}
\noindent\textit{Proof of Lemma \ref{lem2}:}  
Let $\|\cdot\|$ be the  matrix operator norm. By direct calculations can have,
\begin{align*}
E \Big\{ \frac{1}{p} Y \trans B_n(\lambda) Y -\frac{1}{p} \tr B_n(\lambda)\Big \}&=0,\\ 
E \Big\{ \frac{1}{p} Y \trans B_n(\lambda) Y -\frac{1}{p} \tr B_n(\lambda)\Big \}^2&=E \frac{2}{p^2} \tr B^2_n(\lambda) \leq  \frac{2 \|\Sigma\|^2}{ \lambda^2 p} \to 0,\\
E \Big\{ \frac{1}{p} Y \trans B^2_n(\lambda) Y -\frac{1}{p} \tr B^2_n(\lambda)\Big\}&=0,\\ 
E \Big\{ \frac{1}{p} Y\trans  B^2_n(\lambda) Y -\frac{1}{p} \tr B^2_n(\lambda) \Big\}^2&=E \frac{2}{p^2} \tr B^4_n(\lambda) \leq  \frac{2 \|\Sigma\|^4}{ \lambda^4 p} \to 0.
\end{align*}
Here, we used the fact that
\begin{align*}
B_n(\lambda)=\Sigma^{\frac{1}{2}} (S_n+\lambda I)^{-1} \Sigma^{\frac{1}{2}}\ind (\frac{1}{n-2}W_n+ \lambda \Sigma^{-1})^{-1},
\end{align*}
where $W_n \sim \mbox{Wishart}(I_p,n-2)$ is semidefinite and hence 
\begin{align*}
\|B_n(\lambda)\| \leq \frac{1}{\lambda_{min}(\lambda \Sigma^{-1})}=\frac{\|\Sigma\|}{\lambda}.
\end{align*}
Consequently we have,
\begin{align*}
\frac{1}{p} Y \trans B_n(\lambda) Y=\frac{1}{p} \tr B_n(\lambda)+o_p(1),~\frac{1}{p} Y \trans B^2_n(\lambda) Y=\frac{1}{p} \tr B^2_n(\lambda)+o_p(1). 
\end{align*}
Together with Theorem 2 of \cite{wang2015shrinkage} we have, 
\begin{align}\label{lem42}
& \frac{1}{p} \tr B_n(\lambda) \inp R_1(\lambda),~\frac{1}{p} \tr B^2_n(\lambda) \inp R_2(\lambda).
\end{align}
The proof is completed.  \hfill$\fbox{}$

\begin{lem} \label{lem3}
Under the conditions (C1)-(C3),	
	\begin{align*}
	&\mu\trans B_n(\lambda) \mu \inp \Delta^2 H_1(\lambda),~\mu \trans  B^2_n(\lambda)\mu \inp  \Delta^2 H_2(\lambda).
	\end{align*}
\end{lem}
\noindent\textit{Proof of Lemma \ref{lem3}:}  
By \cite{karoui2011geometric} we have,
\begin{align*}
& \mu \trans B_n(\lambda) \mu -\mu \trans \Sigma^{\frac{1}{2}} (\gamma_n \Sigma+\lambda I_p)^{-1} \Sigma^{\frac{1}{2}} \mu \inp 0,\\
&  \mu \trans B^2_n(\lambda) \mu-(1+\epsilon_n)\mu \trans \Sigma^{\frac{1}{2}} (\gamma_n \Sigma+\lambda I_p)^{-1} \Sigma (\gamma_n \Sigma+\lambda I_p)^{-1} \Sigma^{\frac{1}{2}}\mu \inp 0,
\end{align*}
where 
\begin{align*}
\gamma_n=\frac{1}{1+\frac{1}{n-2}  \tr B_n(\lambda)}, ~\epsilon_n= \frac{\gamma_n^2}{n-2}  \tr B^2_n(\lambda).
\end{align*}
From \eqref{lem42}
we have
\begin{align*}
& \gamma_n \inp \frac{1}{1+y R_1(\lambda)}=1-y(1-\lambda m_0(-\lambda))\defby \gamma, \\
& \epsilon_n \inp \frac{y R_2(\lambda)}{(1+y R_1(\lambda))^2}= \frac{y(1-\lambda m_0(-\lambda))}{1-y(1-\lambda
	m_0(-\lambda))}-\frac{y(\lambda m_0(-\lambda)-\lambda^2
	m'_0(-\lambda))}{(1-y(1-\lambda m_0(-\lambda)))^2}\defby \epsilon.
\end{align*}
Consequently it can be routinely shown that
\begin{align*}
& \mu \trans B_n(\lambda) \mu -\mu \trans \Sigma^{\frac{1}{2}} (\gamma \Sigma+\lambda I_p)^{-1}\Sigma^{\frac{1}{2}} \mu \inp 0,\\
&  \mu \trans B^2_n(\lambda) \mu-\frac{1+\epsilon}{\gamma^2} \Big\{ \mu \trans\Big(I_p+\frac{\lambda}{\gamma} \Sigma^{-1}\Big)^{-2} \mu\Big\}  \inp 0,
\end{align*}
By Condition (C3),
\begin{align*}
\mu \trans \Sigma^{\frac{1}{2}} (\gamma \Sigma+\lambda I_p)^{-1}\Sigma^{\frac{1}{2}} \mu=\gamma^{-1} \mu \trans\Big(I_p+\frac{\lambda}{\gamma} \Sigma^{-1}\Big)^{-1} \mu \to\gamma^{-1} h_1\Big(\frac{\lambda}{\gamma}\Big) \Delta^2,
\end{align*}
and 
\begin{align*}
\frac{1+\epsilon}{\gamma^2} \Big\{\mu \trans\Big(I_p+\frac{\lambda}{\gamma} \Sigma^{-1}\Big)^{-2} \mu\Big\}\to\frac{1+\epsilon}{\gamma^2} h_2\Big(\frac{\lambda}{\gamma}\Big) \Delta^2.
\end{align*}

The proof is completed.  \hfill$\fbox{}$

\subsection{Proof of Theorem \ref{thm1}}
By the properties of Gaussian distributions,
\begin{align*}
\bar{X}_1\ind \frac{1}{\sqrt{n_1}} \Sigma^{\frac{1}{2}}Y_1+\mu_1,~\bar{X}_2 \ind \frac{1}{\sqrt{n_2}} \Sigma^{\frac{1}{2}}Y_2+\mu_2,
\end{align*}
where $Y_1,Y_2 \sim N(0,I_p)$ and $Y_1,Y_2$ are independent. We then have,
\begin{align*}
2 \Big(\mu_1-\frac{\bar{X}_1+\bar{X}_2}{2}\Big) \trans \Sigma^{-1}(\bar{X}_1-\bar{X}_2)&=\Delta^2-\frac{2}{\sqrt{n_2}} \mu \trans Y_2+\frac{1}{n_2} Y_2 \trans Y_2-\frac{1}{n_1} Y_1 \trans Y_1   \\
&\inp \Delta^2+y_2-y_1,
\end{align*}
and
\begin{align*}
-2\Big(\mu_2-\frac{\bar{X}_1+\bar{X}_2}{2}\Big)\trans\Sigma^{-1} (\bar{X}_1-\bar{X}_2)&=\Delta^2+\frac{2}{\sqrt{n_1}} \mu \trans Y_1+\frac{1}{n_1} Y_1 \trans Y_1-\frac{1}{n_2} Y_2 \trans Y_2   \\
&\inp \Delta^2+y_1-y_2,
\end{align*}
where $\mu$ is defined as in \eqref{Bn}. 
For the denominator,
\begin{align*}
&(\bar{X}_1-\bar{X}_2)\trans \Sigma^{-1}(\bar{X}_1-\bar{X}_2)\\
=& \Delta^2+\frac{1}{n_1} Y\trans_1 Y_1+\frac{1}{n_2} Y\trans_2 Y_2+2 \mu \trans \Big(\frac{1}{\sqrt{n_1}} Y_1-\frac{1}{\sqrt{n}_2} Y_2\Big)-\frac{2}{\sqrt{n_1 n_2}} Y\trans_1 Y_2\\
\inp& \Delta^2+y_1+y_2. 
\end{align*}
By the Slutsky's theorem and the continuous mapping theorem, the proof is completed. \hfill$\fbox{}$
\subsection{Proof of Theorem \ref{thm2}}
Write 
\begin{align*}
S_n \ind \frac{1}{n-2}\Sigma^{\frac{1}{2}} W_n \Sigma^{\frac{1}{2}},
\end{align*}
where $W_n \sim \mbox{Wishart}(I_p,n-2)$. Then
\begin{align*}
\r_{\text{LDA2}}\ind \Phi\Big(\frac{-\mu \trans W_n^{-1} \mu }{2\sqrt{\mu \trans W_n^{-2} \mu}}\Big).
\end{align*}
By Lemma \ref{lem1}, when $p/n \to y \in (0,1)$, we have
\begin{align*}
& n E \mu \trans W_n^{-1} \mu=\frac{n\Delta^2 }{n-p-3} \to \frac{\Delta^2}{1-y},\\
& n^2 E \mu \trans W_n^{-2} \mu= \frac{n^2(n-3)\Delta^2}{(n-p-2)(n-p-3)(n-p-5)} \to \frac{\Delta^2}{(1-y)^3},
\end{align*}
and 
\begin{align*}
& Var(n \mu \trans W_n^{-1} \mu)=\Delta^4\Big\{ \frac{n^2}{(n-p-3)(n-p-5)}- \frac{n^2}{(n-p-3)^2}\Big\} \to 0,\\
& Var(n^2 \mu \trans W_n^{-2} \mu)=\Delta^4\Big\{ \frac{n^6}{(n-p)^6}(1+o(1))-  \frac{n^4(n-3)^2}{(n-p-2)^2(n-p-3)^2(n-p-5)^2}\Big \}\to 0.
\end{align*}
Therefore, 
\begin{align*}
n \mu \trans W_n^{-1} \mu \inp  \frac{\Delta^2}{1-y},~n^2 \mu \trans W_n^{-2} \mu \inp  \frac{\Delta^2}{(1-y)^3}.
\end{align*}
The proof is completed. \hfill$\fbox{}$
\subsection{Proof of Theorem \ref{thm3}}
We have
\begin{align*}
\bar{X}_1\ind \frac{1}{\sqrt{n_1}} \Sigma^{\frac{1}{2}}Y_1+\mu_1,~\bar{X}_2 \ind \frac{1}{\sqrt{n_2}} \Sigma^{\frac{1}{2}}Y_2+\mu_2, S_n \ind \frac{1}{n-2}\Sigma^{\frac{1}{2}} W_n \Sigma^{\frac{1}{2}},
\end{align*}
where $Y_1,Y_2 \sim N(0,I_p),~W_n \sim \mbox{Wishart}(I_p,n-2)$ and $Y_1,Y_2, W_n$ are independent. Then

\begin{align*}
&-\Big(\mu_1-\frac{\bar{X}_1+\bar{X}_2}{2}\Big)\trans S_n^{-1}(\bar{X}_1-\bar{X}_2)\\
\ind & -\frac{n-2}{2}\Big(\mu-\frac{1}{\sqrt{n_1}}Y_1-\frac{1}{\sqrt{n_2}}Y_2\Big) \trans W_n^{-1}\Big(\mu+\frac{1}{\sqrt{n_1}}Y_1-\frac{1}{\sqrt{n_2}}Y_2\Big)\\
=& -\frac{n-2}{2}\Big(\mu-\frac{1}{\sqrt{n_2}}Y_2\Big) \trans W_n^{-1}\Big(\mu-\frac{1}{\sqrt{n_2}}Y_2\Big)+\frac{n-2}{2n_1} Y_1 \trans W_n^{-1}  Y_1.
\end{align*}
By Lemma \ref{lem1}, we have,
\begin{align*}
E Y_1 \trans W_n^{-1}  Y_1=\frac{1}{n-p-3} E Y_1 \trans Y_1=\frac{p}{n-p-3} \to \frac{y}{1-y},
\end{align*}
and 
\begin{align*}
Var(Y_1 \trans W_n^{-1}  Y_1)=&E (Y_1 \trans W_n^{-1}  Y_1)^2-( E Y_1 \trans W_n^{-1}  Y_1)^2\\
=&\frac{E (Y_1 \trans Y_1)^2}{(n-p-3)(n-p-5)}-\frac{p^2}{(n-p-3)^2}\\
=& \frac{p^2+2p}{(n-p-3)(n-p-5)} -\frac{p^2}{(n-p-3)^2} \to 0.
\end{align*}
Similarly, 
\begin{align*}
& E (n-2)\Big(\mu-\frac{1}{\sqrt{n_2}}Y_2\Big) \trans W_n^{-1}\Big(\mu-\frac{1}{\sqrt{n_2}}Y_2\Big) \\
=& \frac{n-2}{n-p-3} E \Big(\mu-\frac{1}{\sqrt{n_2}}Y_2\Big) \trans \Big(\mu-\frac{1}{\sqrt{n_2}}Y_2\Big)\\
=& \frac{n-2}{n-p-3}\Big(\Delta^2+\frac{p}{n_2}\Big) \to \frac{1}{1-y}(\Delta^2+y_2),
\end{align*}
and 
\begin{align*} 
&Var\Big\{ (n-2)\Big(\mu-\frac{1}{\sqrt{n_2}}Y_2\Big) \trans W_n^{-1}\Big(\mu-\frac{1}{\sqrt{n_2}}Y_2\Big) \Big\}\\
=&\frac{(n-2)^2E \{ \big(\mu-\frac{1}{\sqrt{n_2}}Y_2\big) \trans \big(\mu-\frac{1}{\sqrt{n_2}} Y_2\big)\}^2}{(n-p-3)(n-p-5)}-\frac{(n-2)^2\big(\Delta^2+\frac{p}{n_2}\big)^2}{(n-p-3)^2}\\
=& \frac{(n-2)^2\big(\Delta^4+\frac{(2p+4)\Delta^2}{n_2}  +\frac{p^2+2p}{n_2^2}\big)}{(n-p-3)(n-p-5)}-\frac{(n-2)^2\big(\Delta^2+\frac{p}{n_2}\big)^2}{(n-p-3)^2} \to 0.
\end{align*}
Thus,
\begin{align}
&-\Big(\mu_1-\frac{\bar{X}_1+\bar{X}_2}{2}\Big)\trans S_n^{-1}(\bar{X}_1-\bar{X}_2) \nonumber \\
\ind & -\frac{n-2}{2}\Big(\mu-\frac{1}{\sqrt{n_2}}Y_2\Big) \trans W_n^{-1}\Big(\mu-\frac{1}{\sqrt{n_2}}Y_2\Big)+\frac{n-2}{2n_1} Y_1 \trans W_n^{-1}  Y_1 \nonumber \\
\inp & -\frac{1}{2(1-y)}(\Delta^2+y_2)+\frac{y_1}{2y}  \frac{y}{1-y}=-\frac{\Delta^2-y_1+y_2}{2(1-y)}, \label{pthm31}
\end{align}
and by similar arguments, we have
\begin{align}
&\Big(\mu_2-\frac{\bar{X}_1+\bar{X}_2}{2}\Big)\trans S_n^{-1}(\bar{X}_1-\bar{X}_2) \nonumber \\
\ind & -\frac{n-2}{2}\Big(\mu+\frac{1}{\sqrt{n_1}}Y_1\Big) \trans W_n^{-1}\Big(\mu+\frac{1}{\sqrt{n_1}}Y_1\Big)+\frac{n-2}{2n_2} Y_2 \trans W_n^{-1}  Y_2 \nonumber \\
\inp & -\frac{\Delta^2+y_1-y_2}{2(1-y)}. \label{pthm32}
\end{align}
For the denominator, we have
\begin{align*}
& (\bar{X}_1-\bar{X}_2) \trans S_n^{-1} \Sigma S_n^{-1}(\bar{X}_1-\bar{X}_2)\\
\ind & (n-2)^2\Big(\mu+\frac{1}{\sqrt{n_1}}Y_1-\frac{1}{\sqrt{n_2}}Y_2\Big) \trans W_n^{-2}\Big(\mu+\frac{1}{\sqrt{n_1}}Y_1-\frac{1}{\sqrt{n_2}}Y_2\Big).
\end{align*}
By Lemma \ref{lem1}, we have,
\begin{align*}
& E (n-2)^2\Big(\mu+\frac{1}{\sqrt{n_1}}Y_1-\frac{1}{\sqrt{n_2}}Y_2\Big) \trans W_n^{-2}\Big(\mu+\frac{1}{\sqrt{n_1}}Y_1-\frac{1}{\sqrt{n_2}}Y_2\Big)\\
=& \frac{(n-2)^2(n-3) E \big(\mu+\frac{1}{\sqrt{n_1}}Y_1-\frac{1}{\sqrt{n_2}}Y_2\big) \trans \big(\mu+\frac{1}{\sqrt{n_1}}Y_1-\frac{1}{\sqrt{n_2}}Y_2\big)}{(n-p-2)(n-p-3)(n-p-5)} \\
=& \frac{(n-2)^2(n-3)\big(\Delta^2+\frac{p}{n_1}+\frac{p}{n_2}\big)}{(n-p-2)(n-p-3)(n-p-5)}   \to \frac{\Delta^2+y_1+y_2}{(1-y)^3},
\end{align*}
and
\begin{align*}
& E \Big\{(n-2)^2\Big(\mu+\frac{1}{\sqrt{n_1}}Y_1-\frac{1}{\sqrt{n_2}}Y_2\Big) \trans W_n^{-2}\Big(\mu+\frac{1}{\sqrt{n_1}}Y_1-\frac{1}{\sqrt{n_2}}Y_2\Big)\Big\}^2\\
=&\Big\{\frac{n^6}{(n-p)^6}+o(1)\Big\} E\Big\{ \Big(\mu+\frac{1}{\sqrt{n_1}}Y_1-\frac{1}{\sqrt{n_2}}Y_2\Big) \trans \Big(\mu+\frac{1}{\sqrt{n_1}}Y_1-\frac{1}{\sqrt{n_2}}Y_2\Big)\Big\}^2\\
=&\Big\{\frac{n^6}{(n-p)^6}+o(1)\Big\} \Big\{\Big(\Delta^2+\frac{p}{n_1}+\frac{p}{n_2}\Big)^2+o(1)\Big\}  \to \frac{(\Delta^2+y_1+y_2)^2}{(1-y)^6}.
\end{align*}
Thus, 
\begin{align} \label{pthm33}
(\bar{X}_1-\bar{X}_2) \trans S_n^{-1} \Sigma S_n^{-1}(\bar{X}_1-\bar{X}_2) \inp  \frac{\Delta^2+y_1+y_2}{(1-y)^3}.
\end{align}
Combing \eqref{pthm31}, \eqref{pthm32} and \eqref{pthm33}, the proof is completed. \hfill$\fbox{}$

{\subsection{Proof of  Proposition \ref{prop0}}
With some simple calculation it can be shown that the misclassification rate of $\delta_{\rm LDA}^c(X)$ is }
\begin{align*}
\r_{\text{LDA}}^c=\frac{1}{2}\sum_{j=1}^2 \Phi\left( \frac{ (-1)^j\big[(\mu_j-\frac{\bar{X}_1+\bar{X}_2}{2})\trans S_n^{-1}(\bar{X}_1-\bar{X}_2)
+\alpha\big]}
{\sqrt{(\bar{X}_1-\bar{X}_2) \trans S_n^{-1} \Sigma S_n^{-1}(\bar{X}_1-\bar{X}_2)}}\right).
\end{align*}
where $\alpha=\frac{n-2}{n-p-3}\big(\frac{p}{2 n_1}-\frac{p}{2 n_2}\big)$. Note that from \eqref{pthm31} and the fact that $\alpha\rightarrow \frac{y_1-y_2}{2(1-y)}$, we have
\begin{eqnarray*}
-\bigg[\Big(\mu_1-\frac{\bar{X}_1+\bar{X}_2}{2}\Big)\trans S_n^{-1}(\bar{X}_1-\bar{X}_2)
 +\alpha\bigg]&\inp& -\frac{\Delta^2-y_1+y_2}{2(1-y)}-\frac{y_1-y_2}{2(1-y)}\\
 &=&-\frac{\Delta^2}{2(1-y)}.
\end{eqnarray*}
Similarly, by \eqref{pthm32} we have
\begin{eqnarray*}
	\bigg[\Big(\mu_2-\frac{\bar{X}_1+\bar{X}_2}{2}\Big)\trans S_n^{-1}(\bar{X}_1-\bar{X}_2)
	+\alpha\bigg]&\inp&  -\frac{\Delta^2}{2(1-y)}.
\end{eqnarray*}
Combining with \eqref{pthm33} we have
\begin{eqnarray*}
\r_{\text{LDA}}^c&\inp& \Phi\left(-\frac{\Delta^2}{2(1-y)} \frac{(1-y)^{3/2}}{\sqrt{\Delta^2+y_1+y_2}} \right)\\
&=&\Phi\left(- \frac{\Delta^2}{2\sqrt{\Delta^2+y_1+y_2}} \sqrt{1-y}\right).
\end{eqnarray*}

\subsection{Proof of  Theorem \ref{thm4}}
Write 
\begin{align*}
&T_{1n}=(2\mu_1-\bar{X}_1-\bar{X}_2)\trans (S_n+\lambda I_p)^{-1}(\bar{X}_1-\bar{X}_2),\\
&T_{2n}=-(2\mu_2-\bar{X}_1-\bar{X}_2)\trans (S_n+\lambda I_p)^{-1}(\bar{X}_1-\bar{X}_2),\\
&T_{3n}=(\bar{X}_1-\bar{X}_2)\trans (S_n+\lambda I_p)^{-1}\Sigma (S_n+\lambda I_p)^{-1} (\bar{X}_1-\bar{X}_2).
\end{align*} 
Since 
\begin{align*}
\bar{X}_1\ind \frac{1}{\sqrt{n_1}} \Sigma^{\frac{1}{2}}Y_1+\mu_1,~\bar{X}_2 \ind \frac{1}{\sqrt{n_2}} \Sigma^{\frac{1}{2}}Y_2+\mu_2,
\end{align*}
where $Y_1,Y_2 \sim N(0,I_p)$ and $Y_1,Y_2,S_n$ are independent, we have
\begin{align*}
T_{1n} \ind & \mu \trans B_n(\lambda) \mu-\frac{2}{\sqrt{n_2}} \mu \trans B_n(\lambda) Y_2+\frac{1}{n_2} Y_2 \trans B_n(\lambda) Y_2-\frac{1}{n_1} Y_1 \trans B_n(\lambda)Y_1,\\
T_{2n} \ind & \mu \trans B_n(\lambda) \mu+\frac{2}{\sqrt{n_1}} \mu \trans B_n(\lambda) Y_1+\frac{1}{n_1} Y_1 \trans B_n(\lambda) Y_1-\frac{1}{n_2} Y_2 \trans B_n(\lambda)Y_2,\\
T_{3n} \ind &\Big(\mu+\frac{1}{\sqrt{n_1}} Y_1-\frac{1}{\sqrt{n_2}} Y_2\Big) \trans B^2_n(\lambda) \Big(\mu+\frac{1}{\sqrt{n_1}} Y_1-\frac{1}{\sqrt{n_2}} Y_2\Big)\\
\ind & \mu \trans B^2_n(\lambda) \mu+2\sqrt{\frac{1}{n_1}+\frac{1}{n_2}} \mu \trans B^2_n(\lambda) Y_1+\Big(\frac{1}{n_1}+\frac{1}{n_2}\Big)Y_1\trans B^2_n(\lambda) Y_1.
\end{align*} 
For the leading terms, by Lemmas \ref{lem2} and \ref{lem3},
\begin{align*}
&\frac{1}{n_j} Y_j \trans B_n(\lambda)Y_j \inp y_j R_1(\lambda),~j=1,2\\
&\Big(\frac{1}{n_1}+\frac{1}{n_2}\Big)Y_1\trans B^2_n(\lambda) Y_1 \inp (y_1+y_2) R_2(\lambda)\\
&\mu \trans B_n(\lambda) \mu \inp H_1 (\lambda) \Delta^2,~\mu \trans B^2_n(\lambda) \mu \inp H_2 (\lambda) \Delta^2.
\end{align*}
For the cross sectional terms, 
\begin{align*}
E \Big\{\frac{1}{\sqrt{n_j}} \mu \trans B_n(\lambda) Y_j\Big\} ^2=\frac{1}{n_j} E \mu \trans B^2_n(\lambda) \mu \leq \frac{ \Delta^2\|\Sigma\|^2 }{n_j \lambda^2} \to 0,
\end{align*} 
and similarly,
\begin{align*}
E \Big\{\frac{1}{\sqrt{n}} \mu \trans B^2_n(\lambda) Y_1\Big\} ^2=\frac{1}{n} E \mu \trans B^4_n(\lambda) \mu \leq \frac{ \Delta^2\|\Sigma\|^4 }{n_j \lambda^4} \to 0.
\end{align*}
Thus,
\begin{align*}
\frac{1}{\sqrt{n_1}} \mu \trans B_n(\lambda) Y_1,\frac{1}{\sqrt{n_2}} \mu \trans B_n(\lambda) Y_2,\sqrt{\frac{1}{n_1}+\frac{1}{n_2}} \mu \trans B^2_n(\lambda) Y_1 \inp 0.
\end{align*}
Above all, 
\begin{align*}
& T_{1n} \inp H_1(\lambda) \Delta^2+(y_2-y_1) R_1(\lambda),~T_{2n} \inp H_1(\lambda)+(y_1-y_2) R_1(\lambda),\\
&T_{3n} \inp H_2(\lambda) \Delta^2 +(y_1+y_2) R_2(\lambda).
\end{align*}
By the continuous mapping theorem, the proof is completed. \hfill$\fbox{}$
\subsection{Proof of  Proposition \ref{prop2}}
When
	\begin{align*}
h_1(t)=\int \frac{x}{x+t} dH(x),~h_2(t)=\int \frac{x}{(x+t)^2} dH(x),	
\end{align*}
by the definitions of $H_1(\lambda)$ and $H_2(\lambda)$, we have
\begin{align*}
H_1(\lambda)=& (1+yR_1(\lambda))h_1(\lambda (1+yR_1(\lambda)))=\int \frac{1}{\frac{\lambda}{x}+\frac{1}{1+yR_1(\lambda)}}dH(x),
\end{align*}
and
\begin{align*}
H_2(\lambda)=&\{(1+yR_1(\lambda))^2+y R_2(\lambda)\} \int \frac{x^2}{(x+\lambda (1+yR_1(\lambda))^2)} dH(x)\\
=&\int \frac{1+\frac{y R_2(\lambda)}{(1+yR_1(\lambda))^2}}{(\frac{\lambda}{x}+\frac{1}{1+yR_1(\lambda)})^2}dH(x).
\end{align*}
Compared with Lemma 2 of \cite{wang2015shrinkage},  $H_1(\lambda)=R_1(\lambda)$ and $H_2(\lambda)=R_2(\lambda)$. The convergence of misclassification rate is a direct conclusion of Theorem \ref{thm4}.  The proof is completed. \hfill$\fbox{}$ 


\bibliographystyle{imsart-nameyear}
\bibliography{ref}
\end{document}